\documentclass[11pt,a4paper,twoside]{article}
\usepackage{color}
\usepackage{bm, amsmath, amssymb, amsthm} 

\def\calf{{\cal F}}
\def\<{\langle}
\def\>{\rangle}
\def\eps{\varepsilon}

\def\RR{\mathbb{R}}
\newcommand\codim{\operatorname{codim}}
\newcommand\const{\operatorname{const}}

\newcommand\tr{\operatorname{Tr}}
\newcommand\Div{\operatorname{div}}
\newcommand\id{\operatorname{id}}

\def\Ric{\operatorname{Ric}}
\def\vol{\operatorname{vol}}

\def\eq{\hspace*{-1.5mm}&=&\hspace*{-1.5mm}}
\def\plus{\hspace*{-1.5mm}&+&\hspace*{-1.5mm}}
\def\minus{\hspace*{-1.5mm}&-&\hspace*{-1.5mm}}

\newtheorem{corollary}{Corollary}
\newtheorem{definition}{Definition}
\newtheorem{example}{Example}
\newtheorem{remark}{Remark}
\newtheorem{lem}{Lemma}
\newtheorem{prop}{Proposition}
\newtheorem{thm}{Theorem}

\author{Vladimir Rovenski\footnote{Faculty of Natural Sciences, Mathematical Department,
       University of Haifa, Mount Carmel, 31905 Haifa,  Israel.
       \newline e-mail: {\tt vrovenski@univ.haifa.ac.il}
       }
        \ and \
        Pawe\l \ Walczak\footnote{Faculty of Mathematics and Computer Science, Department of Geometry,
        University of \L\'{o}dz, ul. 22 Banach Str., 90-238  \L\'{o}d\'{z}, Poland.
        e-mail: {\tt pawelwal@math.uni.lodz.pl}}
}

\title{Integral formulae\\ for codimension-one foliated Randers spaces}

\begin{document}

\date{}

\maketitle


\begin{abstract}
Integral formulae
for foliated Riemannian manifolds
provide obstructions for existence of foliations or compact leaves of them with given geometric properties.
This paper continues our recent study and presents new integral formulae and their applications for codimension-one foliated  Randers spaces. The goal is a generalization of Reeb's formula (that the~total mean curvature of the leaves is zero)
and its companion (that twice total second mean curvature of the leaves equals to the total Ricci curvature in the normal direction). We~also extend results by Brito, Langevin and Rosenberg (that total mean curvatures of arbitrary order for a codimension-one foliated  Riemannian manifold of constant curvature don't depend on a foliation).
All of that is done by a comparison of extrinsic and intrinsic curvatures
of the two Riemannian structures which arise in a natural way from a given Randers structure.



\end{abstract}

\section*{Introduction}
\label{sec:intro}

Two recent decades brought increasing interest in Finsler spaces $(M,F)$, especially,
in extrinsic geometry of their hypersurfaces, see \cite{cs2,sh1,sh2}.
Randers metrics $F=\alpha+\beta$, where $\alpha$ is the norm
of a Riemannian structure and $\beta$ a 1-form of $\alpha$-norm smaller than~$1$ on~$M$
(introduced in \cite{ra} and appeared in a solution \cite{brs} of Zermelo's control problem)
are of particular interest, see~\cite{cs}.
Extrinsic geometry of foliated Riemannian manifolds also became popular
since some time (see \cite{rw1} and the bibliography therein).
Among other topics of interest, one can find
so called {\it integral formulae} (i.e., integral relations for invariants of the shape operator of leaves,
e.g. the mean curvatures $\sigma_k\ (1\le k\le m)$, and Riemann curvature, see surveys in \cite{rw1,arw2014}).
Such formulae provide obstructions for existence of foliations or compact leaves of them with given geometric properties.
 The~first known integral formula (by G.\,Reeb, \cite{re}) for codimension-1 foliated closed manifolds tells~us~that
the~total mean curvature $H=\sigma_1$ of the leaves is zero (thus, either $H\equiv0$ or $H(x)H(y)<0$  for some points $x,y\in M$).
Its~counterpart in the case of the second mean curvature is the
(according to our knowledge, obtained for the first time in \cite{no}) formula
\begin{equation}\label{eq61b-init}
 \int_M (2\,\sigma_2-\Ric_{N})\,{\rm d}V_g=0,
\end{equation}
where $N$ is a unit normal to the leaves.
Such formulae were used in \cite{lw} to prove that codimension-one foliations of a closed Riemannian manifold of
either negative Ricci curvature or constant nonzero curvature are far (in a sense defined there) from being totally umbilical, and in \cite{bw} to estimate the energy of a vector field.
An infinite series of integral formulae was provided in~\cite{rw2}: they include the Reeb's formula and generalize
Brito-Langevin-Rosenberg formulae \cite{blr}, which show that total mean curvatures (of arbitrary order $k$) for codimension-one foliations
on a closed $(m+1)$-dimensional manifold of constant sectional curvature $K$ depend only
on $K$, $k$, $m$ and the volume of the manifold, not on a foliation.
In \cite{rw3}, we studied integral formulae for a codimension-one foliation $\calf$ of a closed Finsler space $(M,F)$;
using a unit vector field $\nu$ orthogonal (in the Finsler sense) to the leaves
we defined a new Riemannian structure $g$ on $M$ and derived its Riemann curvature and the shape operator of $\calf$
in terms of $F$.
Using our approach in \cite{rw2}, we produced the integral formulae for $(M, F)$ and for Randers space
$(M,\alpha+\beta)$ with $\beta^\sharp$ (i.e., the $\alpha$-dual of $\beta$) tangent to the leaves.

This paper presents new integral formulae for a codimension-one foliated Randers space.
Section~\ref{sec:prelim} surveys necessary facts and recent results.
Section~\ref{sec:randers} contains our main results,
which generalize the Reeb's formula and \eqref{eq61b-init},
and extend certain formulae in \cite{blr};
in particular, we generalize some results of~\cite{rw3}.
All integral formulae of this paper hold when the foliation
and the 1-form, both  are defined outside a finite union of closed submanifolds
of codimension $\ge 2$ under convergence of some integrals  (as in Lemma~\ref{L-LuW-2} in what follows),
leaving details to the readers.
The singular case is important since there exist plenty of manifolds which admit no (smooth) codimension-one foliations,
while all of them admit such foliations (and non-singular 1-forms $\beta$) outside some ``set of singularities".

\section{Preliminaries}
\label{sec:prelim}

We work with a closed manifold $M$ equipped with a codimension-one foliation defined
on $M\setminus\Sigma$, where $\Sigma$ is a (possibly empty) union of pairwise disjoint closed submanifolds
$\Sigma_i$ of variable codimensions $\ge2$.
Briefly, we say that our foliation admits singularities at points of $\Sigma$.
For Randers spaces (with metrics $F=\alpha+\beta$) we assume also that $\beta$
admits singularities, i.e., is defined on $M\setminus\Sigma$.

\begin{lem}[see Lemma~2 in \cite{lw2}]\label{L-LuW-2}
Let $\Sigma_1$,  $\codim \Sigma_1\ge2$, be a closed submanifold of a Riemannian manifold $(M,a)$,
and $X$ a vector field on $M\setminus\Sigma_1$ such that $\int_M\|X\|^2\,{\rm d}V_a<\infty$.
Then $\int_M\Div X\,{\rm d}V_a=0$.
\end{lem}

For $\sigma_k\ (k\ge2)$ the singular case is also considered in \cite[Theorem~2]{rw2}.

\subsection{The Minkowski and Randers norms}
\label{sec:M-R-norms}

\begin{definition}[see \cite{sh2}]\rm
A \textit{Minkowski norm} on a vector space $V^{m+1}$ is a function $F:V^{m+1}\to[0,\infty)$ with the following properties (of regularity, positive 1-homogeneity and strong convexity):

M$_1:$ $F\in C^\infty(V^{m+1}\setminus \{0\})$,\qquad
M$_2:$ $F(\lambda\,y)=\lambda F(y)$ for all $\lambda>0$ and $y\in V^{m+1}$,

M$_3:$ For any $y\in V^{m+1}\setminus \{0\}$, the following symmetric bilinear form is positive definite on $V^{m+1}:$
\begin{equation}\label{E-acsiom-M3}
 g_y(u,v)=\frac12\,\frac{\partial^2}{\partial s\,\partial t}\,\big[F^2(y+su+tv)\big]_{|\,s=t=0}\,.
\end{equation}
\end{definition}

\noindent
By (M$_2$) and (M$_3$), $g_{\lambda y}=g_{y}\ (\lambda>0)$, and
$\{y\in V^{m+1}: F(y) \le 1\}$ is a strictly convex set.
Note that
\begin{equation}\label{E-g-F2}
 g_y(y,v)=\frac12\frac{\partial}{\partial t}\,\big[F^2(y+tv)\big]_{|\,t=0},\quad
 g_y(y,y)=F^2(y).
\end{equation}
For Minkowski norms, the following symmetric trilinear form $C$ is called the \textit{Cartan torsion}:
\begin{equation}\label{E-Ctorsion}
  C_y(u,v,w)=\frac12\,\frac{\partial}{\partial t}\,\big[g_{y+tw}(u,v)\big]_{|\,t=0}\quad
 {\rm where} \quad y\in V^{m+1}\setminus\{0\},\ u,v,w\in V^{m+1}\,.
\end{equation}
The homogeneity of $F$ implies
\[
 C_y(u,v,w)=\frac14\,\frac{\partial^3}{\partial r\,\partial s\,\partial t}\,\big[F^2(y+ru+sv+tw)\big]_{|\,r=s=t=0},\quad
 C_{\lambda y}=\lambda^{-1}C_y\quad (\lambda>0).
\]
Moreover, we have $C_y(y,\cdot\,,\cdot\,)=0$. The \textit{mean Cartan torsion} is defined by
 $I_{\,y}(u):=\tr C_{\,y}(\cdot\,,\cdot\,, u)$.
Let $(b_i)$ be a basis for $V^{m+1}$ and $(\theta^i)$ the dual basis in $(V^{m+1})^*$.
The \textit{Busemann-Hausdorff volume form} is defined by
\[
 {\rm d}V_F=\sigma_F(x)\,\theta^1\wedge\dots\wedge\theta^{m+1},\quad {\rm where}\quad \sigma_F=\frac{\vol\mathbb{B}^{m+1}}{\vol B^{m+1}},
\]
where
$\mathbb{B}^{m+1}:=\{y=y^i b_i\in V^{m+1}: \|y\|^2=\sum_i (y^i)^2<1\}$,
and $\vol B^{m+1}$ is the Euclidean volume of the subset
$B^{m+1}:=\{y\in V^{m+1}: F(y^i b_i)<1\}$ of $V$. The~\textit{distortion} of $F$ is defined by
 $\tau(y)=\log((\det g_{ij}(y))^{1/2}/{\sigma_F})$.
It has the $0$-homogeneity property: $\tau(\lambda\,y)=\tau(y)\ (\lambda>0)$,
and $\tau=0$ for Riemannian spaces. The~\textit{angular form} is defined by
 $h_y(u,v)=g_y(u,v) - F(y)^{-2}g_y(y,u)\,g_y(y,v)$.
A~vector $n\in V^{m+1}$ is \textit{normal} to a hyperplane $W\subset V^{m+1}$ if $g_n(n,w)=0\ (w\in W)$.
There are exactly two normal directions to $W$ which are opposite when $F$ is \textit{reversible},
i.e., $F(-y) = F(y)\ (y\in V^{m+1})$.

\begin{definition}[see, for example, \cite{sh2}]\rm
 Let
 $a(\cdot\,,\cdot)= \<\cdot\,,\cdot\>$ be a scalar product and
 $\alpha(y)=\|y\|_\alpha=\sqrt{\<y,y\>}$ for $y\in \RR^{m+1}$ the corresponding Euclidean norm on $\RR^{m+1}$.
 If $\beta$ is a linear form on $\RR^{m+1}$ with the property
 $\|\beta\,\|_\alpha<1$,
 then the following nonnegative function $F$ is called the \textit{Randers norm}:
\begin{equation*}
 F(y) = \alpha(y) + \beta(y) = \sqrt{\<y,y\>}+\beta(y).
\end{equation*}
\end{definition}

For Randers norm
on $\RR^{m+1}$, the bilinear form $g_y$ is positive definite and obeys, see \cite{sh2} again,
\begin{eqnarray}\label{E-F001}
\nonumber
 g_y(u,v)\eq \alpha^{-2}(y)(1+\beta(y))\,\<u,v\> +\beta(u)\,\beta(v) \\
 \minus \alpha^{-3}(y)\,\beta(y)\,\<y,u\>\,\<y,v\> +\alpha^{-1}(y)\,\big(\beta(u)\,\<y,v\> +\beta(v)\,\<y,u\>\big)\,,\\
\label{E-F001b}
 \det g_y \eq (F(y)/\alpha(y))^{m+2}\det a.
\end{eqnarray}
Let $N\in \RR^{m+1}$ be a unit normal to a hyperplane $W$ in $\RR^{m+1}$ with respect to $\<\cdot\,,\cdot\>$, i.e.,
\[
 \<N,w\>=0\quad (w\in W),\qquad
  \alpha(N)=\|N\|_\alpha=\sqrt{\<N,N\>}=1.
\]
Let $n$ be a vector $F$-normal to $W$, i.e.,
 $g_n(n,v) = 0\ (v\in W)$,
lying in the same half-space as $N$ and such that $\|n\|_\alpha=\alpha(n)=1$.
Set
\[
 g(u,v):=g_n(u,v),\quad u,v\in\RR^{m+1}.
\]
Then $g(n,n)=F^2(n)$, see \eqref{E-g-F2}, and  $F(n)=1+\beta(n)$.
For $y=n$ formula \eqref{E-F001} yields
\begin{equation}\label{E-c-value0}
 g(u,v) = (1+\beta(n))\<u,v\> +\beta(u)\,\beta(v) - \beta(n)\,\<n,u\>\,\<n,v\> +\beta(u)\,\<n,v\> +\beta(v)\,\<n,u\>.
\end{equation}
The 'musical isomorphisms' $\sharp$ and $\flat$ will be used
for rank one and symmetric rank 2 tensors
on Riemannian manifolds.
For~example, if $\beta$ is a 1-form on $\RR^{m+1}$ and $v\in \RR^{m+1}$ then
$\<\beta^\sharp,u\>=\beta(u)$ and $v^\flat(u) =\<v,u\>$ for any $u\in \RR^{m+1}$.
The~tangent component of a vector, say $\beta^\sharp$, will be denoted by $\beta^{\sharp\top}$,
its dual 1-form is~$\beta^\top$.

\begin{lem}\label{L-c-value} We have
\begin{eqnarray}\label{E-c-value}
 n \eq \hat c\,N-\beta^\sharp,\quad \mbox{\rm or, equivalently,}\quad n = c\,N-\beta^{\sharp\top}, \\
\label{E-c-value2}
 g(u,v) \eq c\,\hat c\,(\<u,v\> -\beta(u)\,\beta(v)),\quad u,v\in W\,,\\
\label{E-c-value2c}
 g(n,n) \eq (c\,\hat c)^{\,2},
\end{eqnarray}
where $c:=(1-\|\beta^{\sharp\top}\|^2_\alpha)^{\frac12}>0$ and $\hat c = c+\beta(N)$.
 The vector $\nu=(c\,\hat c)^{-1}n$ is an $F$-unit normal to~$W$.
\end{lem}

\proof Assuming $u=n$, from \eqref{E-c-value0} and $g(n,v)=0$ we find
\begin{equation}\label{E-gnv}
 (1+\beta(n))\,\<n + \beta^{\sharp},\,v\> =0.
\end{equation}
Note that
 $|\beta(n)|=|\<\beta^\sharp,n\>|\le\alpha(\beta^\sharp)\,\alpha(n)<1$; hence, $1+\beta(n)>0$.
We find from \eqref{E-gnv} with $v\in W$ that $n+\beta^{\sharp}=\hat c\,N$ for some $\hat c>0$.
Using $1 = \<n,n\> = \hat c^{\,2}-2\,\hat c\,\beta(N) +\|\beta\,\|^2_\alpha$, we get two values
 $\hat c = \beta(N)\pm c$,
from which $\beta(N) + c$ is positive, that proves \eqref{E-c-value}$_1$.
In view of $\beta^\sharp=\beta^{\sharp\top}+\beta(N)N$ this is equivalent to~\eqref{E-c-value}$_2$.
Thus, \eqref{E-c-value2c} follows from $g(n,n)=(1+\beta(n))^2$ and
\[
 1+\beta(n)=1+\beta(\hat c\,N-\beta^\sharp)=c\,\hat c.
\]
 Note that $F(n)=c\,\hat c$. Finally, \eqref{E-c-value2} follows from \eqref{E-c-value0}.
\qed

\begin{lem}\label{L-zZ}
 If $u,U\in W$ and
 $g(u,v)=\<U,v\>$ for all $v\in W$
then
\begin{equation}\label{E-u-U}
 (c\,\hat c)\,u = U+c^{-2}\beta^\top(U)\,\beta^{\sharp\top}.
\end{equation}
\end{lem}

\proof
By \eqref{E-c-value2},
 $g(u,v) = c\,\hat c\,\<u -\beta(u) \beta^{\sharp\top},\,v\>$ holds.
By conditions,
since $u,U$ and $\beta^{\sharp\top}$ belong to $W$, we obtain
 $u -\beta(u)\beta^{\sharp\top} = (c\,\hat c)^{-1} U$.
Applying $\beta$, we obtain $\beta(u)=(c\,\hat c)^{-1}c^{-2}\beta(U)$ and then~\eqref{E-u-U}.
\qed

\subsection{Finsler spaces}

Let $M^{m+1}$ be a connected smooth manifold and $TM$ its tangent bundle.
A~\textit{Finsler structure} $F$ on $M$ is a
family of Minkowski norms in tangent spaces $T_pM$ which depend smoothly on a point $p\in M$.
 The {covariant derivative} of a vector field $u(t)$ along a curve $c(t)$ in $M$ is given by
\[
 D_{\dot c}\,u = \{\dot u^i + \Gamma^i_{kj}(\dot c)\,\dot c^k\,u^j\}\,\partial_{x^i\,|\,c}\,,
\]
where
 $\Gamma^i_{kj}=\frac12\,g^{il}\big(\frac{\partial g_{jl}}{\partial x^k}
 +\frac{\partial g_{kl}}{\partial x^j} -\frac{\partial g_{jk}}{\partial x^l}\big)$
are homogeneous of 0-degree functions on $TM_0:=TM\setminus\{0\}$, and $g_{ij}(y)=\frac12\,[F^2]_{y^i y^j}(y)$,
compare \eqref{E-acsiom-M3}.
A vector field $u$ along a curve $c$ is \textit{parallel} if $D_{\dot c}\,u\equiv0$.
A~curve $c$ is called a {geodesic} if
the tangent vector $u=\dot c$ is parallel along itself: $D_{\dot c}\,\dot c =0$.
A Finsler metric on $M$ is called a \textit{Berwald metric} if in any coordinate system $(x,y)$ in $TM_0$,
the Christoffel symbols $\Gamma^i_{jk}$ are functions on $x\in M$ only;
such Finsler spaces are modeled on a single Minkowski space.
Berwald metrics are characterized among Randers ones $F = \alpha + \beta$,
by the property: $\beta$ is parallel with respect to $\alpha$, see \cite[Theorem~2.4.1]{sh2}.

Let $c_y$ be a geodesic with $\dot c_y(0)=y\in T_pM$.
The~{exponential~map} $\exp_p: y\mapsto c_y(1)$ (by~homogeneity, $\exp_p(ty)=c_y(t)$ for $t>0$)
is smooth on $TM_0$ and $C^1$ at the origin with $d(\exp_p)_{|\,0}=\id_{\,T_pM}$, see \cite{sh1}.
A $C^\infty$ map ${\mathcal H}:(-\eps,\eps)\times[0,1]\to M$ is called
a {geodesic variation} of a geodesic $c(t),\ 0\le t\le 1$, if ${\mathcal H}(0,t)=c(t)$ and for each $s\in(-\eps,\eps)$,
the curve $c_s(t):={\mathcal H}(s,t)$ is a geodesic.
The variation field $Y(t):=\frac{\partial{\mathcal H}}{\partial s}(0,t)$ along $c$ obeys the {Jacobi~equation}:
\begin{equation}\label{E-RC-Jacobi}
 D_{\dot c}D_{\dot c}\,Y +R_{\,\dot c}(Y) =0
\end{equation}
for some $y$-dependent $g_y$-self-adjoint (1,1)-tensor $R_y$, called the \textit{Riemann curvature in a direction} $y\in T_pM\setminus\{0\}$.
 By \eqref{E-RC-Jacobi}, $R_y(y)=0$ and $R_{\,\lambda y}=\lambda^2 R_y\ (\lambda>0)$.
Let $\{e_i\}_{1\le i\le m+1}$ be a $g_y$-orthonormal basis for $T_pM$ such that $e_{m+1}=y/F(y)$,
and let $P_i={\rm span}\{e_i,y\}$ for some $y\in T_pM$.
The \textit{Ricci curvature} is a (positive homogeneous of degree 2) function on $TM_0$:
\[
 \Ric_{\,y}=\sum\nolimits_{\,i=1}^m g_y(R_y(e_i),e_i).
 \]

\subsection{Codimension-one foliated Finsler spaces}
\label{sec:intform}

Given a transversally oriented codimension-one foliation $\calf$ of $(M^{m+1},F)$,
 there exists a globally defined $F$-normal (to the leaves) smooth vector field $n$
which defines a Riemannian metric $g:=g_n$ with the Levi-Civita connection $\nabla$.
Then $g(n,u)=0\ (u\in T\calf)$ and $g(n,n)=F^2(n)$, see~\eqref{E-c-value2c},
and $\nu=n/F(n)$ is a $F$-unit normal.
The shape operator $A^g:T\calf\to T\calf$ of $\calf$ with respect to the metric $g$ is given
by
\begin{equation}\label{E-Ag}
 A^g(u)=-{\nabla_u\,\nu}\quad (u\in T\calf).
\end{equation}
Let $L$ be the leaf through a point $p\in M$, and $\rho$ the local distance function to $L$
in a neighborhood of~$p$.
Denote by $\hat\nabla$ the Levi-Civita connection of the (local again) Riemannian metric $\hat g:= g_{\,\nabla\rho}$. Note that $\nabla\rho=\nu$ on~$L$.
The~{shape operator} $A:T\calf\to T\calf$ (self-adjoint for $g$) is defined by
\[
 A(u)=-{\hat\nabla_u\,\nu}\quad (u\in T\calf).
\]
Let $C_\nu^\sharp$ be a $(1,1)$-tensor $g$-dual to the symmetric bilinear form $C_\nu(\cdot\,,\cdot\,,\nabla_\nu\,\nu)$.
Note that $C_n^\sharp=\hat c^{\,3}C_\nu^\sharp$\,.

In \cite{rw3}, we applied the variational approach to express the Riemann curvature of $g$
in terms of Riemann curvature and the Cartan torsion of $F$.

\begin{thm}[see Theorem 3.4 in \cite{rw3}]\label{T-unit-n}
Let $\nu$ be an $F$-unit normal to a codimension-one foliation of a Finsler space $(M^{m+1},F)$.
The Riemann curvatures in the $\nu$-direction of $F$ and $g$, the shape operators and volume forms are related~by
\begin{eqnarray}\label{E-dt-Rm-Rm}
\nonumber
 &&\hskip-7mm g((R_{\,\nu} - R^g_{\,\nu})(u),v) = -C_{\nu}\big(A^g(u)+\frac12\,C_\nu^\sharp(u), v, \nabla_{\nu}\,\nu\big)
 +2\,(\nabla_{\nu}C_{\nu})(u,v,\nabla_{\nu}\,\nu)\\
 &&\hskip4mm -\,C_{\nu}\big(u, A^g(v) +\frac12\,C_\nu^\sharp(v), \nabla_{\nu}\,\nu\big)
 +C_{\nu}\big(u, v, \nabla^2_{\nu,\nu}\,\nu -C_\nu^\sharp(\nabla_{\nu}\,\nu)\big) ,\\
\label{E-dt-A2}
 &&\hskip0mm A - A^g = C_\nu^\sharp,\qquad {\rm d}V_g = e^{\tau(\nu)}\,{\rm d}V_F.
\end{eqnarray}
\end{thm}

In case of a Riemannian foliation (i.e., when the vector field $\nu$ is geodesic: $\nabla_\nu\,\nu=0$) from \eqref{E-dt-Rm-Rm} we obtain $R_{\,\nu} = R^g_{\,\nu}$, see also \cite[Proposition~6.2.2]{sh2}.

Invariants $\sigma_{\lambda} (A_1, \ldots, A_k)$ of a set of real $m\times m$ matrices
are discussed briefly in Section~\ref{sec:inv}.
They generalize elementary symmetric functions $\sigma_i(A)$ of a single symmetric matrix~$A$ (i.e., $k=1$).
Recall that $\sigma_1(A)=\tr A$.

\begin{thm}[see Theorem 3.6 in \cite{rw3}]\label{thm:main}
If $\mathcal F$ is a codimension-one foliation with a unit normal $\nu$ on a closed
$F$-locally symmetric Finsler manifold $(M^{m+1}, F)$, then for any $0\le k\le m$ one has
\begin{equation}\label{eq:main}
 \int_M \sum\nolimits_{\,\|\lambda\|=k}\sigma_{\lambda}\left(B_1, \ldots B_k\right)\,{\rm d}V_F=0,
\end{equation}
where
 $B_{2k}=\frac{(-1)^k}{(2k)!}\,(R_{\nu})^k,\
 B_{2k+1}=\frac{(-1)^k}{(2k+1)!}\,(R_{\nu})^k A$.
\end{thm}

The formulae (\ref{eq:main}) for initial values of $k$, $k=1,2$, read as follows:
\begin{eqnarray}\label{eq61}
 \int_M \sigma_1(A)\,{\rm d}V_F = 0,\qquad
 \int_M \big( \sigma_2(A) - \frac12\Ric_{\,\nu}\big)\,{\rm d}V_F = 0.
\end{eqnarray}
Recall (see the Introduction) that \eqref{eq61} are known for foliated Riemannian spaces.

Next corollary of Theorem~\ref{thm:main} generalizes result for Riemann manifolds in~\cite{blr}.

\begin{corollary}[see Corollary 3.9 in \cite{rw3}]\label{thm:003}
Let ${\mathcal F}$ be a transversally oriented codimension-one foliation on a
closed
Finsler mani\-fold $(M^{m+1},F)$
with a unit normal $\nu$ and condition $R_{\nu} = K\,I_m$. Then, for any $0\le k\le m$,
\begin{equation}\label{eq5f-b}
 \int_M \sigma_k(A)\,{\rm d}V_F = \left\{
 \begin{array}{cc}
 K^{k/2} \genfrac{(}{)}{0pt}{1}{\,m/2\,}{k/2}
 \,{\rm Vol}_F(M), & m,\,k \ {\rm even}, \\
 0, & m \ {\rm or}\ k \ {\rm odd}.
 \end{array}\right.
\end{equation}
\end{corollary}

\begin{remark}\label{sec:nonsym}\rm
Theorem~\ref{thm:main} and Corollary~\ref{thm:003} are valid for
a
foliation with singularities of codimension $\ge k$, due to Theorem~2 and Corollary~4 of~\cite{rw2}.
Moreover, the compactness of $M$ can be replaced by weaker conditions that
$M$~has finite $F$-volume, and `bounded geometry' in the following~sense:
\begin{equation}\label{E-bounded}
 \sup\nolimits_{\,M} \|R_{\,\nu}\|_F<\infty,\quad \sup\nolimits_{\,M} \|A\|_F<\infty.
\end{equation}
\end{remark}

\section{Codimension-one foliated Randers spaces}\label{sec:randers}

This section generalizes results in \cite{rw3}, where the case of $\beta(N)=0$ has been studied.
As before, write $\langle\cdot,\cdot\rangle$ -- a  Riemannian metric on~$M^{m+1}$.
 Let $\calf$ be a transversally oriented codimension-one foliation of a Randers space $(M^{m+1},F)$:
\[
 F(y)=\sqrt{\<y,y\>} +\beta(y),\quad  \|\beta\,\|_\alpha<1,\quad \beta^\sharp\in\Gamma(TM).
\]
Let $N$ be a unit $\alpha$-normal vector field to $\calf$,
and $n$ an $F$-normal vector field to $\calf$ with the property $\<n,n\>=1$.
Let $\bar\nabla$ be the Levi-Civita connection of $\langle\cdot,\cdot\rangle$,
and $\nabla$ the Levi-Civita connection of the metric $g=g_n$ on~$M$.
By \cite[(1.15) \& (1.19)]{cs},
\begin{eqnarray}
\nonumber
 \tau(n) \eq \frac12\,(m+2)\log{\frac{1+\beta(n)}{1-\|\beta\|_\alpha^2}}
 =\frac{m+2}2\,\log{\frac{c}{2\,c-\hat c}}\,,\\
 \label{E-Iy}
 I_n(u) \eq \tr C_{\,n}(\cdot\,,\cdot\,, u) =\frac{m+2}{2\,c\,\hat c}\,\<\,\beta^\sharp - (c\,\hat c -1)\,u,\ n\>\,,
\end{eqnarray}
 where $c=\sqrt{1-\|\beta^{\sharp\top}\|^2_\alpha}>0$ and $\hat c=c+\beta(N)>0$, see Lemma~\ref{L-c-value}.
Observe that
\begin{equation*}
 C_n(u,v,w) = \frac1{m+2}\,\big( I_n(u)\,h_n(v,w)+I_n(v)\,h_n(u,w)+I_n(w)\,h_n(u,v)\big)\,,
\end{equation*}
where
$h_n(u,v) = c\,\hat c\,( \<u,v\> -\<u,n\>\,\<v,n\>)$
is the angular form, see \cite[(1.11) \& (1.20)]{cs}.
 We have $\sigma_F=(1-\|\beta^{\sharp}\|^{2}_\alpha)^{\frac{m+2}2}\sqrt{\det a_{ij}}$, see \cite
{cs},
and
$\sqrt{\det g_{ij}(n)}=(c\,\hat c)^{\frac{m+2}2}\sqrt{\det a_{ij}}$, see \eqref{E-F001b}.
Thus, the canonical volume forms of metrics $g$ and $\langle\cdot,\cdot\rangle$ satisfy
\begin{equation}\label{E-F001vol}
 {\rm d}V_F=
 (1-\|\beta^{\sharp}\|^{2}_\alpha)^{\frac{m+2}2}{\rm d}V_a,\ \
 {\rm d} V_g = (c\,\hat c)^{\frac{m+2}2}{\rm d} V_a,\ \
 {\rm d}V_F= (1-\|\beta^{\sharp}\|^{2}_\alpha)^{\frac{m+2}2}{\rm d}V_g.
\end{equation}
Recall that $\nu =(c\,\hat c)^{-1} n$.
Let $Z=\nabla_\nu\,\nu$ and $\bar Z=\bar\nabla_N\,N$ be the curvature vectors of $\nu$- and $N$- curves for
$g$ and $\langle\cdot,\cdot\rangle$, respectively. In the case of $\beta^{\sharp\top}\ne0$,
let $X^{\bot\beta}$ be the projection of $X\in \Gamma(T\calf)$ on~$\beta^{\sharp\bot}$:
\begin{equation}\label{E-X-beta-bot}
 X^{\bot\beta}=X-\<X,\ \beta^{\sharp\top}\>\,\|\beta^{\sharp\top}\|^{-2}_\alpha\,\beta^{\sharp\top}\,.
\end{equation}
Notation \eqref{E-X-beta-bot} will be used in decompositions of matrices $\tilde B = B +\sum_i B_i$,
where $B_i$ are rank 1 matrices of the form $U^{\bot\beta}\otimes\beta^{\sharp\top}$,
$(U^{\bot\beta})^\flat\otimes\beta^{\sharp\top}$ and $f\cdot\beta^{\top}\otimes\beta^{\sharp\top}$
for some $U\in T\calf$. The invariants of $\tilde B$ and $B$ are close in the sense, see Appendix.

\subsection{The shape operators of $g$ and $\langle\cdot,\cdot\rangle$}

The derivative $\bar\nabla u:TM\to TM$ and its conjugate $(\bar\nabla u)^t: TM\to TM$
are $(1,1)$-tensors defined by $(\bar\nabla u)\,(v)=\bar\nabla_v\,u$
and $\<(\bar\nabla\,u)^t(v),w\>=\<v,(\bar\nabla\,u)(w)\>$ for $v,w\in TM$.
The \textit{deformation tensor},
 $2\,\overline{\rm Def}_u
 =\bar\nabla u+(\bar\nabla u)^t$,
measures the degree to which the flow of a vector field $u$ distorts the metric $\langle\cdot,\cdot\rangle$.
The~same notation $\overline{\rm Def}_u$ will be used for its $\langle\cdot,\cdot\rangle$-dual $(1,1)$-tensor.
Set $\overline{\rm Def}^\top_u(v)=(\overline{\rm Def}_u(v))^\top$.

\begin{prop}\label{L-Dx}
The shape operators of $g$ and $\langle\cdot,\cdot\rangle$ satisfy on $\calf$ the following:
\begin{eqnarray}\label{E-A-bar-A-initial}
\nonumber
 && c\,A^g = \bar A -\frac12\,c^{-1}\hat c^{\,-2}(\hat c\,N-\beta^{\sharp})(c\,\hat c)\,I_m
 +\hat c^{\,-1}\,(\overline{\rm Def}_{\beta^\sharp})_{\,|T\calf}^\top
 +\,\frac12\,\big( U -\bar A(\beta^{\sharp\top})\,\big)\otimes\beta^\top \\
 &&\hskip-2mm
 +\,\frac12\,c^{-2}\Big( \bar A(\beta^{\sharp\top}) -\<\bar  A(\beta^{\sharp\top}),\,\beta^{\sharp\top}\>\,\beta^{\sharp\top}
 +\,2\,\hat c^{\,-1}(\overline{\rm Def}_{\beta^\sharp}\,\beta^{\sharp\top})^\top
 +U +\beta(U)\,\beta^{\sharp\top}\Big) \,\!^\flat\otimes\beta^{\sharp\top},
\end{eqnarray}
where $\,U=\hat c^{\,-1}(\bar\nabla_{\hat c\,N-\beta^{\sharp}}\,\beta^{\sharp\top})^\top -c\bar Z$\,.
At points $p\in M$ with $\beta^{\sharp\top}(p)\ne0$ we~get
\begin{eqnarray}\label{E-A-bar-A}
\nonumber
 && c\,A^g = \bar A -\frac12\,c^{-1}\hat c^{\,-2}(\hat c\,N-\beta^{\sharp})(c\,\hat c)\,I_m
 +\hat c^{\,-1}\,(\overline{\rm Def}_{\beta^\sharp})_{\,|T\calf}^\top
 +\,\frac12\,\big( U -\bar A(\beta^{\sharp\top})\,\big)\,^{\bot\beta}\otimes\beta^\top \\
 &&
 +\frac12\,c^{-2}\big(\,2\,\hat c^{\,-1}(\overline{\rm Def}_{\beta^\sharp}\,\beta^{\sharp\top})^\top
  +( U +\bar A(\beta^{\sharp\top}))^{\bot\beta}\big)\,\!^\flat\otimes\beta^{\sharp\top}
  +\frac1{c^2(1-c^2)}\,\beta(U)
 \,\beta^\top\!\otimes\beta^{\sharp\top}.
\end{eqnarray}
\end{prop}

\proof
By well-known formula for the Levi-Civita connection of $g$
and use of the equalities $g(u,n)=0=g(v,n)$ and $g([u,v],n)=0$ we have
\begin{equation}\label{E-LC-g}
 2\,g(\nabla_u\,n, v) = n(g(u,v)) +g([u,n],v) +g([v,n],u)\quad (u,v\in T\calf).
\end{equation}
Assume $\bar\nabla_X^\top\,u=\bar\nabla_X^\top\,v=0$ for all $X\in T_pM$ at a given point $p\in M$.
Using \eqref{E-c-value0}--\eqref{E-c-value2}, we obtain
\begin{eqnarray*}
 n(g(u,v)) \eq n\big(c\,\hat c\,(\<u,v\> -\beta(u)\,\beta(v))\big)
  = n(c\,\hat c)(\<u,v\>-\beta(u)\beta(v))\\
 \minus c\,\hat c\,\big(\beta(u)(\bar\nabla_n(\beta^\top))(v) +(\bar\nabla_n(\beta^\top))(u)\beta(v)\big),\\
 g([u,n],v) \eq c\,\hat c\,\big(\<\,[u,n],v\> +\<[u,n],\,n\,\>\beta(v)\big)\\
 \eq -c\,\hat c\,\<\hat c\,\bar A(u) +\bar\nabla_u\,\beta^\sharp,\,v\>
 +c\,\hat c^{\,2}\<\bar A(\beta^{\sharp\top}) +c\bar Z,\, u\>\,\beta(v),\\
 g([v,n],u) \eq c\,\hat c\,\big(\<\,[v,n],u\> +\beta(u)\<[v,n],\,n\,\>\big)\\
 \eq - c\,\hat c\,\<\hat c\,\bar A(v) +\bar\nabla_v\,\beta^\sharp,\,u\>
 +c\,\hat c^{\,2}\beta(u)\,\<\bar A(\beta^{\sharp\top}) +c\bar Z,\, v\> .
\end{eqnarray*}
Substituting the above into \eqref{E-LC-g}, we find
\begin{eqnarray}\label{E-DDD}
\nonumber
 2\,g(\nabla_u\,n, v) \eq n(c\,\hat c)(\<u,v\>-\beta(u)\beta(v)) -2\,c\,\hat c^{\,2}\<\bar A(u), v\>
 -2\,c\,\hat c\,\<\overline{\rm Def}_{\beta^\sharp}(u),v\> \\
\nonumber
 \minus c\,\hat c\,\big(\beta(u)(\bar\nabla_n(\beta^\top))(v) +(\bar\nabla_n(\beta^\top))(u)\beta(v)\big)\\
 \plus c\,\hat c^{\,2}\big(\beta(v)\,\<\bar A(\beta^{\sharp\top}) +c\bar Z,\, u\>
 +\beta(u)\,\<\bar A(\beta^{\sharp\top}) +c\bar Z,\, v\>\big).
\end{eqnarray}
Assume $g(\nabla_u\,n, v)=\<\mathfrak{D}(u),\,v\>$, where $\mathfrak{D}: T\calf\to T\calf$ is a linear operator.
Using Lemma~\ref{L-zZ} and $g(\nabla_u\,n,v) =-c\,\hat c\,g(A^g(u),v)$, see \eqref{E-Ag},
we get from \eqref{E-DDD} the following:
\begin{equation}\label{E-Dx}
 -2\,(c\,\hat c)^{\,2} A^g(u) = 2\,\mathfrak{D}(u)+c^{-2}\<2\,\mathfrak{D}(u),\,\beta^{\sharp\top}\>\,\beta^{\sharp\top},
\end{equation}
where
\begin{eqnarray}\label{E-Dx2}
\nonumber
 2\,\mathfrak{D}(u) \eq n(c\,\hat c)(u-\beta(u)\beta^{\sharp\top}) -2\,c\,\hat c^{\,2}\bar A(u)
 -2\,c\,\hat c\,(\overline{\rm Def}_{\beta^\sharp}(u))^\top
 -c\,\hat c\,\big(\beta(u)(\bar\nabla_n\,\beta^{\sharp\top})^\top\\
 \plus (\bar\nabla_n(\beta^\top))(u)\,\beta^{\sharp\top}\big)
 +c\,\hat c^{\,2}\big(\<\bar A(\beta^{\sharp\top}) +c\bar Z,\, u\>\,\beta^{\sharp\top}
 +\beta(u)\,(\bar A(\beta^{\sharp\top}) +c\bar Z)\big).
\end{eqnarray}
In particular, using $\<\beta^{\sharp\top},\beta^{\sharp\top}\>=1-c^2$ we get
\begin{eqnarray}\label{E-Dx3}
\nonumber
 &&\hskip-10mm \<2\,\mathfrak{D}(u),\beta^{\sharp\top}\> = n(c\,\hat c)\,c^2\beta(u) -2\,c\,\hat c^{\,2}\<\bar A(\beta^{\sharp\top}),\,u\>
 -2\,c\,\hat c\,\<\,\overline{\rm Def}_{\beta^\sharp}(\beta^{\sharp\top}),\,u\> \\
\nonumber
 && 
 -\,c\,\hat c\,\big( \beta(u)\<\bar\nabla_n(\beta^{\sharp\top}),\,\beta^{\sharp\top}\>
 +(1-c^2)(\bar\nabla_n(\beta^\top))(u)\big)\\
 && +\,c\,\hat c^{\,2}\big((1-c^2)\<\bar A(\beta^{\sharp\top}) +c\bar Z,\, u\>
    +\beta(u)\,\<\bar A(\beta^{\sharp\top})+c\bar Z,\,\beta^{\sharp\top}\> \big).
\end{eqnarray}
 From \eqref{E-Dx}, \eqref{E-Dx2} and \eqref{E-Dx3} we obtain
\begin{eqnarray*}
 c\,A^g(u) \eq \bar A(u) -\frac12\,c^{-1}\hat c^{\,-2}n(c\,\hat c)(u-\beta(u)\beta^{\sharp\top})
 +\hat c^{\,-1}(\overline{\rm Def}_{\beta^\sharp}(u))^\top\\
\nonumber
 && +\,\frac12\,\hat c^{\,-1}\big(\beta(u)(\bar\nabla_n\,\beta^{\sharp\top})^\top
 +(\bar\nabla_n(\beta^\top))(u)\,\beta^{\sharp\top}\big)\\
 && -\,\frac12\,\big(\<\bar A(\beta^{\sharp\top}) +c\bar Z,\,u\>\,\beta^{\sharp\top}
 +\beta(u)\,(\bar A(\beta^{\sharp\top}) +c\bar Z)\big)\\
 && +\,\frac12\,c^{-2}\Big( 2\,\hat c^{\,-1}\<\,\overline{\rm Def}_{\beta^\sharp}(\beta^{\sharp\top}),\,u\>
 +2\,\<\bar A(\beta^{\sharp\top}),\,u\> -c\,\hat c^{\,-2}n(c\,\hat c)\,\beta(u)\\
\nonumber
 && +\,\hat c^{\,-1}\big(\beta(u)\<\bar\nabla_n(\beta^{\sharp\top}),\,\beta^{\sharp\top}\>
 +(1-c^2)(\bar\nabla_n(\beta^\top))(u)\big)\\
 && -\,\big((1-c^2)\<\bar A(\beta^{\sharp\top}) +c\bar Z,\, u\>
 +\beta(u)\,\<\bar A(\beta^{\sharp\top})+c\bar Z,\,\beta^{\sharp\top}\> \big) \Big)\beta^{\sharp\top}.
\end{eqnarray*}
Reducing terms with factors $1-c^2$ and $n(c\,\hat c)\beta(u)$, we obtain  \eqref{E-A-bar-A-initial}.
For $\,\beta^{\sharp\top}\ne0$ we
apply \eqref{E-X-beta-bot} with $X=\bar A(\beta^{\sharp\top})$ and $X=U$, and find \eqref{E-A-bar-A} using
\begin{eqnarray*}
 && U\otimes\beta^{\top} +c^{-2}\,(U+\beta(U)\,\beta^{\sharp\top})\,^\flat\otimes\beta^{\sharp\top} \\
 && =U^{\bot\beta}\otimes\beta^{\top} +c^{-2}\,(U^{\bot\beta})^\flat\otimes\beta^{\sharp\top}
 +\frac{2\,\beta(U)}{c^2(1-c^2)}\,\beta^\top\otimes\beta^{\sharp\top}.\quad\qed
\end{eqnarray*}

\begin{example}\label{C-shape}\rm
Let $\beta(N)=0$ $($i.e., $\hat c=c<1)$ on $M$. Then \eqref{E-A-bar-A} reads as
\begin{eqnarray}\label{E-A-bar-Aold}
\nonumber
 &&\hskip-5mm c\,A^g = \bar A -c^{-2}(c\,N-\beta^\sharp)(c)\,I_m
 +c^{-1}\,(\overline{\rm Def}_{\beta^\sharp})_{\,|T\calf}^\top
 +\frac12\,\big( U -\bar A(\beta^{\sharp})\,\big)\,^{\bot\beta}\otimes\beta\\
 &&\hskip-7mm +\,\frac12\,c^{-2}\big(\,2\,c^{-1}(\overline{\rm Def}_{\beta^\sharp}\,\beta^{\sharp})^\top
  +( U +\bar A(\beta^{\sharp}))^{\bot\beta}\big)\,\!^\flat\otimes\beta^{\sharp}
 +\frac{\beta(U)}{c^2(1-c^2)}\,\beta\otimes\beta^{\sharp},
\end{eqnarray}
where
 $U=(\bar\nabla_{N-c^{-1}\beta^{\sharp}}\,\beta^{\sharp})^\top -c\bar Z$
 and
 $\beta(U)=-(cN-\beta^\sharp)(c)-c\,\beta(\bar Z)$.
This coincides with {\rm \cite[Proposition~3]{rw3}}.
Moreover,
\eqref{E-A-bar-Aold} for $\bar\nabla\beta=0$ reads as
\begin{equation*}
 c\,A^g = \bar A -\frac12\,\big(\bar A(\beta^\sharp)+c\,\bar Z\,\big)\,^{\bot\beta}\otimes\beta
 +\frac12\,c^{-2}\big(\bar A(\beta^\sharp)^\flat -c\bar Z^\flat\big)\,^{\bot\beta}\otimes\beta^\sharp
 -\frac{c^{-1}\beta(\bar Z)}{1-c^2}\,\beta\otimes\beta^\sharp.
\end{equation*}
\end{example}

\subsection{The Riemann curvature}

In this section, we find a relationship between Riemann curvature of metrics $g$ and $\langle\cdot,\cdot\rangle$ on a~Randers space.
For $\beta(N)=0$, the results of this section have been obtained in \cite{rw3}.

\begin{prop}\label{P-ZZ} We have
\begin{equation}\label{E-ZZ}
 Z = (c\,\hat c)^{-1}\bar Z -c^{-1}\hat c^{\,-2}\,\bar\nabla^\top\hat c
 +c^{-3}\hat c^{\,-1}\beta(\bar Z-\hat c^{\,-1}\,\bar\nabla^\top\hat c)\,\beta^{\sharp\top}
\end{equation}
and
\begin{equation}\label{E-Csharp}
 (c\,\hat c)C^\sharp_n = \bar C +c^{-2}(\beta^\top\circ\bar C)\otimes\beta^{\sharp\top},
\end{equation}
where
\begin{eqnarray*}
 &&\hskip-6mm 2\,\bar C = \beta^\top\otimes\bar Z + \beta^{\sharp\top}\otimes\bar Z^\flat
 -\hat c^{\,-1} (\beta^\top\otimes\bar\nabla^\top c + \beta^{\sharp\top}\otimes(\bar\nabla^\top c)^\flat) \\
 &&\hskip-6mm +(c\,\hat c)^{-1}\big(
 (\hat c -2\,c^{-1})\,\beta^{\sharp\top}(\hat c) +\,(c-\hat c^{\,-1})\,n(\hat c)
 +\beta(\bar Z)( c\,\hat c -\hat c^{\,2} +2\,c^{-1}\hat c -1)\big)I_m\\
 &&\hskip-6mm +(c\,\hat c)^{\,-1}\big(
 (2\,c^{-1}-3\,\hat c)\,\beta^{\sharp\top}(\hat c)+(\hat c^{\,-1}-3\,c)\,n(\hat c)
 +\beta(\bar Z)(3\,\hat c^{\,2} -3\,c\,\hat c -2\,c^{-1}\hat c+1)\big)\,\beta^\top\otimes\,\beta^{\sharp\top}.
\end{eqnarray*}
\end{prop}

\proof
Extend $X\in T_p\calf$ at a point $p\in M$ onto a neighborhood of $p$ with the property $(\bar\nabla_Y\,X)^\top=0$ for any $Y\in T_pM$. By the well known formula for the Levi-Civita connection, we obtain at $p$:
 $g(Z,X) = g([X,\nu],\nu)$.
Then, using
$\nu=\hat c^{\,-1}N-(c\,\hat c)^{-1}\beta^{\sharp\top}$ and $[X, fY]=X(f)Y+f[X, Y]$ we get
\begin{eqnarray*}
 g([X,\nu],\nu) \eq \hat c^{\,-3} X(\hat c)\,\big(c^{-1} g(N,\beta^{\sharp\top}) -g(N,N)\big)
 +\hat c^{\,-2} \big( g([X,N],\,N) -c^{-1} g([X,N],\beta^{\sharp\top})\big).
\end{eqnarray*}
Note that
\[
 [X,N] = \bar\nabla_X N -\bar\nabla_N X =-\bar A(X) -\<\bar\nabla_N X,\,N\>\,N
 = -\bar A(X) +\<\bar Z,\,X\>\,N
\]
and $N=\hat c\,\nu +c^{-1}\beta^{\sharp\top}$. Then, by Lemma~\ref{L-c-value} and the equalities
\begin{eqnarray*}
 g(\beta^{\sharp\top},\beta^{\sharp\top}) \eq c^2(\<\beta^{\sharp\top},\,\beta^{\sharp\top}\> -\beta(\beta^{\sharp\top})^2) = c^3\hat c\,(1-c^2),\\
 g(N,\beta^{\sharp\top}) \eq g(\hat c\,\nu+c^{-1}\beta^{\sharp\top},\,\beta^{\sharp\top})
 = c^{-1}g(\beta^{\sharp\top},\beta^{\sharp\top}) = c^2\,\hat c\,(1-c^2),\\
 g(N,N) \eq g(\hat c\,\nu+c^{-1}\beta^{\sharp\top},\,\hat c\,\nu+c^{-1}\beta^{\sharp\top}) =
 \hat c^{\,2} +c^{-2}g(\beta^{\sharp\top},\beta^{\sharp\top}) = \hat c^{\,2} + c\,\hat c\,(1-c^2),
\end{eqnarray*}
we obtain
\begin{eqnarray*}
 g([X,N],\,N) \eq \<\bar Z, X\>\,g(N,N) -c^{-1}\<\bar A(\beta^{\sharp\top}),\,X\>
 =\hat c\,\<(\hat c +c(1-c^2))\bar Z -c^2\bar A(\beta^{\sharp\top}),\ X\>,\\
 g([X,N], \beta^{\sharp\top}) \eq \<\bar Z, X\>\,g(N,\beta^{\sharp\top}) -\<\bar A(\beta^{\sharp\top}),\,X\> =
 c^2\hat c\,\<(1-c^2)\bar Z -c\bar A(\beta^{\sharp\top}),\ X\>.
\end{eqnarray*}
Hence,
 $g(Z,X) = \<\bar Z -\hat c^{\,-1}\,\bar\nabla\,\hat c,\ X\>$.
Applying Lemma~\ref{L-zZ}, we get \eqref{E-ZZ}.
Using definition of $I_n$ and $h_n$, \eqref{E-ZZ} and a bit of help from Maple program we find
\begin{eqnarray*}
 &&\hskip-6mm 2\,C_n(u,v, Z)=\beta(v)\<u,\bar Z\>+\beta(u)\<v,\bar Z\>-\hat c^{\,-1}(\beta(v)\,u(\hat c)
 +\beta(u)\,v(\hat c))\\
 &&\hskip-6mm +\,(c\,\hat c)^{-1}\big((\hat c -2\,c^{-1})\,\beta^{\sharp\top}(\hat c) +(c-\hat c^{\,-1})\,n(\hat c)\big)\<u,v\> \\
 &&\hskip-6mm +\,(c\,\hat c)^{\,-1}\big((2\,c^{-1}-3\,\hat c)\,\beta^{\sharp\top}(\hat c)+(\hat c^{\,-1}-3\,c)\,n(\hat c)\big)\beta(v)\beta(u)\\
 &&\hskip-6mm +\,(c\,\hat c)^{-1}\big(( c\,\hat c {-}\hat c^{\,2} {+}2\,c^{-1}\hat c -1) \<u,v\>
 {+}(3\,\hat c^{\,2} -3\,c\,\hat c -2\,c^{-1}\hat c +1)\,\beta(v)\beta(u)\big)\beta(\bar Z).
\end{eqnarray*}
We have $g(C^\sharp_n(u),v)=\<\bar C(u),v\>$, where $C^\sharp_n$ is $g$-dual to $C_n(\cdot,\cdot, \nabla_n\,n)$, and
\begin{eqnarray*}
 &&\hskip-6mm 2\,\bar C(u) = \<u,\bar Z\>\,\beta^{\sharp\top} +\beta(u)\bar Z
 -\hat c^{\,-1}(u(\hat c)\,\beta^{\sharp\top} +\beta(u)\,\bar\nabla^\top\hat c)\\
 &&\hskip-6mm +\,(c\,\hat c)^{-1}\big((\hat c -2\,c^{-1})\,\beta^{\sharp\top}(\hat c)+(c-\hat c^{\,-1})\,n(\hat c)\big)\,u \\
 &&\hskip-6mm +\,(c\,\hat c)^{\,-1}\big((2\,c^{-1}-3\,\hat c)\,\beta^{\sharp\top}(\hat c)
 +(\hat c^{\,-1}-3\,c)\,n(\hat c)\big)\,\beta(u)\,\beta^{\sharp\top}\\
 &&\hskip-6mm +\,(c\,\hat c)^{-1}\big(( c\,\hat c -\hat c^{\,2} +2\,c^{-1}\hat c -1)\,u
 +(3\,\hat c^{\,2} -3\,c\,\hat c -2\,c^{-1}\hat c +1)\,\beta(u)\,\beta^{\sharp\top}\big)\beta(\bar Z).
\end{eqnarray*}
Then,  we can apply Lemma~\ref{L-zZ} to get \eqref{E-Csharp}.
\qed

\begin{corollary}\label{C-sharp-n}\rm
(i)~Let $\bar\nabla\beta=0$ and $\beta(N)=\const$, then $\bar Z=0$ provides $C^\sharp_n=0$.

(ii)~Let $m>3$, $\beta(N)=\const\ge0$ and $\|\beta\|_\alpha=\const$, then $C^\sharp_n=0$ if and only if $\bar Z=0$.
\end{corollary}

\proof
(i)~Since $c$ and $\hat c$ are constant, and by Proposition~\ref{P-ZZ}, $\bar C=0$, then $C^\sharp_n=0$.

(ii)~Let $y:=\beta(N)=\hat c-c\in(-1,1)$. The roots $y_1<0<y_2$ of the function
\[
 f_c: y \to c\,\hat c -\hat c^{\,2} +2\,c^{-1}\hat c  -1 = -y^2 - c^{-1}(c^2-2)\,y +1
\]
with parameter $c\in(0,1]$ are $y_{1,2}=\frac1{2\,c}\,(2-c^2\pm\sqrt{4+c^4})$.
Note that $y_2>1$ and $y_1>-1$ for $0<c\le1$.
Hence, $f_c(y)>0$ for $0\le y<1$, while for any $c\in(0,1]$ there exists $\tilde y\in(-1,0)$
such that $f_c(\tilde y)=0$.
In other words, $c\,\hat c -\hat c^{\,2}+2\,c^{-1}\hat c -1 >0$ when $\beta(N)\ge0$.
 By~our assumptions,
\begin{eqnarray*}
 && 2\,\bar C = \beta^\top\otimes\bar Z + \beta^{\sharp\top}\otimes\bar Z^\flat \\
 &&\hskip-4mm +\,(c\,\hat c)^{-1}\beta(\bar Z)\big(( c\,\hat c -\hat c^{\,2}+2\,c^{-1}\hat c -1)\,I_m
 {+}(3\,\hat c^{\,2} -3\,c\,\hat c -2\,c^{-1}\hat c +1)\,\beta^\top{\otimes}\beta^{\sharp\top}\big).
\end{eqnarray*}
Hence, $C^\sharp_n=0$, see \eqref{E-Csharp}, reads
\begin{eqnarray}\label{E-CZ}
\nonumber
 && \beta(\bar Z)( c\,\hat c -\hat c^{\,2} +2\,c^{-1}\hat c -1)\,I_m
 = -c\,\hat c\,(\beta^\top\otimes\bar Z + \beta^{\sharp\top}\otimes\bar Z^\flat) \\
 &&\hskip-4mm -\beta(\bar Z)(3\,\hat c^{\,2} -3\,c\,\hat c -2\,c^{-1}\hat c +1)\,\beta^\top\otimes\beta^{\sharp\top}
 -2\,c^{-1}\hat c\,(\beta^\top\circ\bar C)\otimes\beta^{\sharp\top}.
\end{eqnarray}
If $\beta^{\sharp\top}=0$ then the right hand side of \eqref{E-CZ} vanishes and we get $\beta(\bar Z)=0$; hence,
$C^\sharp_n=0$, see also Remark~\ref{sec:beta0}.
Assume now that $\beta^{\sharp\top}\ne0$.
Since the matrix in the left hand side of \eqref{E-CZ} is conformal,
while the matrix in the right hand side of\eqref{E-CZ} has the form $\omega\otimes\beta^{\sharp\top} -(c\,\hat c)\bar Z^{\,\bot\beta}\otimes\beta^\top$ and rank $\le 3$, for $m>3$ we obtain
\begin{equation}\label{E-2cond}
 \beta(\bar Z)=0,\quad
 \beta^\top\otimes\bar Z + \beta^{\sharp\top}\otimes\bar Z^\flat
 +2\,c^{-2}(\beta^\top\circ\bar C)\otimes\beta^{\sharp\top} =0.
\end{equation}
By \eqref{E-2cond}$_1$, $\bar Z\perp\beta^{\sharp\top}$; thus, \eqref{E-2cond}$_2$ yields
$\bar Z=0$ (that is, $\calf$ is a Riemannian foliation for the metric $\langle\cdot,\cdot\rangle$) and $\bar C=0$.
The converse claim follows
from \eqref{E-Csharp} and the definition of $\bar C$.
\qed

\begin{remark}\label{C-3beta}\rm
For a codimension-one foliation of $(M,a)$ we have \cite{rw3}:
\begin{eqnarray}\label{E-nablaZZ}
 &&\hskip-8mm \<\bar\nabla_u \bar Z,v\> = \<\bar\nabla_v\bar Z,u\>,\qquad
 g(\nabla_u Z,v) = g(\nabla_v Z,u)\qquad (u,v\in T\calf),\\
\label{E-RnRnu-R}
 &&\hskip-8mm \bar R_N = (\overline{\rm Def}_{\bar Z})_{\,|T\calf}^\top
 +\bar\nabla_N \bar A -\bar A^2 -\bar Z^\flat\otimes\bar Z.
\end{eqnarray}
In \cite{cs}, $R_y$ is expressed (using coordinate presentations) through $\bar R_y$ for $y\in TM$.
If $\bar\nabla\beta=0$ $($i.e., $F$ is a Berwald structure$)$ then $R_y = \bar R_y$.
Alternative formulas with relationship between $R_\nu$ and $\bar R_\nu$ follow from \eqref{E-RnRnu-R} and
similar formula for $g$,
where $A^g$ and $Z$ are expressed using $\bar A$ and $\bar Z$ given in Propositions~\ref{L-Dx} and \ref{P-ZZ}.
\end{remark}

Given a transversely oriented codimension-1 foliation $\calf$ of an arbitrary closed Finsler
manifold $(M^{m+1}, F)$, denote by $k_1,k_2,\ldots, k_m\ (k_1\le k_2\le\ldots \le k_m)$ the
principal curvatures (eigenvalues of the shape operator $A$) of the leaves of $\calf$.
If $M$ is oriented and $V_F$ is the Finsler volume form on $M$, then one can consider the integral
 $U^F_\calf = \int_M \sum\nolimits_{\,i<j} (k_i-k_j)^2\,{\rm d}V_F$,
which measures ``how far from umbilicity" is $\calf$
(see also \cite[Example~2.6]{rw1} for Riemannian case).
Similar measure of non-umbilicity (with different powers of $k_i - k_j$ which made it conformally invariant)
for
foliated Riemannian manifolds has been considered in \cite[Section~4.1]{lw}.

\begin{thm}\label{T-umb}
Let $\bar\nabla\beta=0$ on $(M,a)$ and the Randers metric $F=\alpha+\beta$ has
$\Ric_\nu\le -r<0$. Then
\begin{equation}\label{E-int-U-0}
 U^F_\calf\ge (1-\|\beta^\sharp\|^2)^{\frac{m+2}2}m\,r\int_M c^{-2}\,{\rm d}V_a.
\end{equation}
\end{thm}

\proof One may show that
\[
 \sum\nolimits_{\,i<j} (k_i-k_j)^2=m\tr(A^2)-(\tr A)^2=(m-1)\,\sigma^2_1(A) -2\,m\,\sigma_2(A).
\]
Hence, and by integral formula \eqref{eq61}$_2$,
\begin{eqnarray}\label{E-int-U-1}
\nonumber
 && U^F_\calf \ge \int_M [(m-1)\,\sigma^2_1(A) -2\,m\,\sigma_2(A)]\,{\rm d}V_F \ge -m\int_M 2\,\sigma_2(A)\,{\rm d}V_F\\
 &&\hskip6mm = -m\int_M \Ric_\nu{\rm d}V_F.
\end{eqnarray}
By condition $\bar\nabla\beta^\sharp=0$ we have $\|\beta^\sharp\|_\alpha=\const$ and
$\bar R(X,Y)\beta^\sharp=0\ (X,Y\in TM)$. Using equality
\[
 \overline{\Ric}_{\,n}
 =\overline{\Ric}_{\,\hat c\,N-\beta^{\sharp}}=\hat c^{\,2}\,\overline{\Ric}_{\,N} +\overline{\Ric}_{\,\beta^\sharp}
 -2\,\hat c\sum\nolimits_{\,i}\bar R(N,b_i,\beta^\sharp,b_i),
\]
we obtain
 ${\Ric}_{\,\nu}=(c\,\hat c)^{\,-2}\,{\Ric}_{\,n} =(c\,\hat c)^{\,-2}\,\overline\Ric_n =c^{-2}\,\overline{\Ric}_{\,N}$.
 From \eqref{E-int-U-1}, where the volume form is
 ${\rm d}V_F=(1-\|\beta^{\sharp}\|^{2}_\alpha)^{\frac{m+2}2}\,{\rm d}V_a$, see \eqref{E-F001vol},
we find
\begin{equation*}
 U^F_\calf
 \ge -(1-\|\beta^\sharp\|_\alpha^2)^{\frac{m+2}2}m\int_M c^{-2}\,\overline{\Ric}_{\,N}\,{\rm d}V_a,
\end{equation*}
which reduces to \eqref{E-int-U-0} since our assumption $\Ric_\nu\le -r<0$.
\qed

\smallskip

Let $\Sigma$ be a  union of pairwise disjoint closed submanifolds $\Sigma_i\subset M$ of codimensions $\ge 2$.
Following \cite{bw} for Riemannian case,
define the energy of a unit vector field $X$ on $M\setminus\Sigma$ by the formula
\[
 {\mathcal E}(X) = \frac 12\int_M \big(\dim M+\|DX\|_F^2\big)\,{\rm d}V_F
 =\frac{m+1}2\,{\rm Vol}_F(M)+\frac12\int_M \|DX\|_F^2\,{\rm d}V_F.
\]
Let $X=\nu$ be a unit normal to a codimension-one foliation. Using the inequality
 $\|D\,\nu\,\|_F^2\ge \frac{2}{m}\,\sigma_2(A)$,
see \cite{bw} for Riemannian case,
Lemma~\ref{L-LuW-2} and integral formula \eqref{eq61}$_2$, we get the following.

\begin{thm}\label{R-vfield}
Let $(M,\alpha+\beta)$ be a codimension-one foliated Randers space with $\bar\nabla\beta=0$.
Then
\begin{equation}\label{E-int-vfield-1}
 {\mathcal E}(\nu)\ge (1-\|\beta^\sharp\|^2)^{\frac{m+2}2}\Big(\frac{m+1}2\,{\rm Vol}_a(M)
 +\frac1{2m}\int_M c^{-2}\,\overline{\Ric}_{\,N}\,{\rm d}V_a\Big).
\end{equation}
\end{thm}

\begin{remark}\label{R-22}\rm
Recall that generally (i.e., when  $N(\beta)\ne0$), $c^2=1 - \|\beta^\top\|^2\ne\const$ and $1 - \|\beta\|^2$ are not the same quantities in \eqref{E-int-vfield-1}.
If $m\ge2$ then equality holds in \eqref{E-int-vfield-1} if and only if $\nu$ is geodesic and $A=\lambda\,I_m$.
One can get an obvious corollary of \eqref{E-int-vfield-1} when $(M,a)$ is a round sphere:
if $c=\const$ then
\[
 {\mathcal E}(\nu)\ge (1-\|\beta^\sharp\|^2)^{\frac{m+2}2}\,\frac{(m+1)\,c^2+1}{2\,c^2}\,{\rm Vol}_a(S^{m+1}).
\]
One~can also drop the condition $\bar\nabla\beta=0$ (in Theorems~\ref{T-umb} and \ref{R-vfield})
and use the formula
$\Ric_n=\overline\Ric_n +\Theta(n)$ for a certain (explicitly given in \cite[p.~54]{cs}) function $\Theta$ on $TM_0$.
\end{remark}

\subsection{Around the Reeb formula}

Basing on \eqref{eq:main} -- \eqref{eq5f-b}, one may produce a sequence of similar formulae for Randers spaces.
We~will discuss first two of them (i.e., for $\sigma_1$ and $\sigma_2$).
In \cite{re}, G. Reeb proved that the total mean curvature of the leaves of a~codimen\-sion-one foliation
on a closed Riemannian manifold equals zero.
The following formula, see \cite[Lemma~2.5]{rw1}, with any $f\in C^2(M)$,
reduces to the Reeb formula when $f=\const\ne0$:
\begin{equation}\label{E-reeb-f}
 \int_M(f\sigma_1(\bar A)-N(f)) \,{\rm d}V_a =0.
\end{equation}
 Recall that
 $\bar Z=\bar\nabla_NN$  is the curvature
 of $N$-curves for $\langle\cdot,\cdot\rangle$,~and
 $c^2=1-\|\beta^{\sharp\top}\|^2_\alpha,\
 \hat c = c+\beta(N)$.
Results of this section are valid for a closed manifold equipped with a codimension-one foliation
and 1-form with singularities of codimension $\ge 2$, see Lemma~\ref{L-LuW-2}.
Moreover, a closed manifold may be replaced by a complete manifold of finite
volume with bounded geometry, see~\eqref{E-bounded}.

\begin{thm}\label{T-1-1}
Let $(M^{m+1},\,\alpha+\beta)$ be a codimension-one foliated closed Randers space. Then
\begin{equation}\label{E-IF1-Randers-gen0fin}
 \int_M (c\,\hat c)^{\frac{m}2}c^{-2}(\hat c-c) \big(
 c\,N(c) +c\,\beta(\bar Z) +\<\bar A(\beta^{\sharp\top}),\beta^{\sharp}\>\,\big)\,{\rm d}V_a =0\,.
\end{equation}
Moreover, if
$c$ and $\beta(N)\ne0$ are constant
then
\begin{equation}\label{E-IF1-Randers-gen}
 \int_M \<\bar A(\beta^{\sharp\top})+c\,\bar Z,\ \beta^{\sharp}\>\,{\rm d}V_a =0.
\end{equation}
\end{thm}

\proof We calculate
\begin{eqnarray}\label{E-tr-Def}
\nonumber
 \tr\,(\overline{\rm Def}_{\beta^\sharp})_{\,|T\calf}^\top
 \eq\sum\nolimits_{\,i=1}^{\,m}\<\bar\nabla_i\,\beta^\sharp,\,b_i\>
 =\overline{\Div}\,\beta^\sharp - \<\bar\nabla_N\,(\beta^{\sharp\top}\!+\beta(N)N),\,N\>\\
\nonumber
 \eq\overline{\Div}\,\beta^\sharp +\beta(\bar Z) - N(\beta(N)),\\
 \<\overline{\rm Def}_{\beta^\sharp}(\beta^{\sharp\top}),\beta^{\sharp\top}\>\eq
 \<\bar\nabla_{\beta^{\sharp\top}}(\beta^{\sharp\top}+\beta(N)N),\beta^{\sharp\top}\>
 = -c\,\beta^{\sharp\top}(c) -\beta(N)\<\bar A(\beta^{\sharp\top}),\beta^{\sharp}\>.
\end{eqnarray}
Tracing \eqref{E-A-bar-A-initial}, we then obtain
\begin{eqnarray}\label{E-trA-1}
\nonumber
 c\,\sigma_1(A^g)
 \eq \sigma_1(\bar A) -\frac m2\,c^{-1}\hat c^{\,-2}\,(\hat c\,N{-}\beta^{\sharp})(c\,\hat c)
 +\hat c^{\,-1}\big(\,\overline{\Div}\,\beta^\sharp {+}\beta(\bar Z) {-} N(\beta(N)) \big) \\
\nonumber
 \plus\frac12\,(\beta(U)-\<\bar A(\beta^{\sharp\top}),\beta^{\sharp}\>)
 +\frac12\,c^{-2}\big( c^2\<\bar A(\beta^{\sharp\top}),\beta^{\sharp}\> \\
\nonumber
 \minus 2\,\hat c^{\,-1}(c\,\beta^{\sharp\top}(c) +\beta(N)\<\bar A(\beta^{\sharp\top}),\beta^{\sharp}\>)
 +(2-c^2)\beta(U) \big)\\
\nonumber
 \eq \sigma_1(\bar A) -\frac m2\,c^{-1}\hat c^{\,-2}\,(\hat c\,N-\beta^{\sharp})(c\,\hat c)
 +\hat c^{\,-1}\,\overline{\Div}\,\beta^\sharp - (\hat c-c)(c\,\hat c)^{-1}\beta(\bar Z) \\
 \minus (\hat c-c)(c\,\hat c)^{-1}N(c) -\hat c^{\,-1}N(\hat c)
 -c^{-2}\hat c^{\,-1}\beta(N)\<\bar A(\beta^{\sharp\top}),\beta^{\sharp}\> .
\end{eqnarray}
 From \eqref{E-trA-1}, \eqref{eq61}$_1$ for $\langle\cdot,\cdot\rangle$ and $g$,
and using ${\rm d}V_g=(c\,\hat c)^{\frac{m+2}2}\,{\rm d}V_a$, see \eqref{E-F001vol}, we get
\begin{eqnarray*}
\nonumber
 &&\hskip-7mm\int_M (c\,\hat c)^{\frac{m+2}2}c^{-1}\Big(\,
 \sigma_1(\bar A) -\frac m2\,c^{-1}\hat c^{\,-2}\,(\hat c\,N-\beta^{\sharp})(c\,\hat c)
 +\hat c^{\,-1}\,\overline{\Div}\,\beta^\sharp -(\hat c-c)(c\,\hat c)^{-1}\beta(\bar Z) \\
 &&\hskip0mm -\,(\hat c-c)(c\,\hat c)^{-1}N(c) -\hat c^{\,-1}N(\hat c)
 -(\hat c-c)\,c^{-2}\hat c^{\,-1}\<\bar A(\beta^{\sharp\top}),\beta^{\sharp}\>
 \Big)\,{\rm d}V_a = 0.
\end{eqnarray*}
The above, the Divergence Theorem and equality
 $f\,\overline\Div\,\beta^\sharp = \overline\Div\,(f\,\beta^\sharp) -\beta^\sharp(f)$
with $f=(c\,\hat c)^{\frac{m}2}$ yield
\begin{eqnarray}\label{E-IF1-Randers-gen0}
\nonumber
 &&\hskip-5mm \int_M \,(c\,\hat c)^{\frac{m+2}2}c^{-1}
 \Big(\sigma_1(\bar A) -\frac m2\,(c\,\hat c)^{\,-1}\,N(c\,\hat c) -(\hat c-c)(c\,\hat c)^{-1}\beta(\bar Z) \\
 &&\hskip0mm -\,(\hat c-c)(c\,\hat c)^{-1}N(c) -\hat c^{\,-1}N(\hat c)
 -(\hat c-c)\,c^{-2}\hat c^{\,-1}\<\bar A(\beta^{\sharp\top}),\beta^{\sharp}\> \Big)\,{\rm d}V_a = 0,
\end{eqnarray}
which is the Reeb formula when $\beta=0$.
Applying \eqref{E-reeb-f}  we obtain \eqref{E-IF1-Randers-gen0fin}.
Note that for $\beta=0$, we have $c=1=\hat c$; hence, \eqref{E-IF1-Randers-gen0} reduces to the Reeb's formula.
If $c\,\hat c=\const$ then \eqref{E-IF1-Randers-gen0fin} reduces to~\eqref{E-IF1-Randers-gen}.
\qed

\begin{remark}\label{R-b0}\rm
The following application of \eqref{E-IF1-Randers-gen} seems to be interesting.
Let $\bar Z=0$ and a unit vector field $X\in\Gamma(T\calf)$ be an eigenvector of $\bar{A}$ corresponding to an eigenvalue $\lambda:M\to\RR$. By Theorem~\ref{T-1-1},
the vector field $\beta^\sharp=\eps' X+\eps N$, where $\eps =\const\in(-1,1)$ and $\eps'=\const\in(0,\sqrt{1-\eps^2})$,
obeys \eqref{E-IF1-Randers-gen}.
Note that $c^2=1-\eps^2-(\eps')^2$, $\beta(N)=\eps$ and $c\,\hat c=1+\eps$.
Thus, assuming $\eps\ne0$, we~get $\int_M \lambda\,{\rm d}V_a =0$.
Consequently, either $\lambda\equiv0$ on $M$ or $\lambda(x)\cdot\lambda(y)<0$ for some points $x$ and $y$ of $M$.
This implies the classical Reeb formula
$\int_M \sigma_1(\bar A)\,{\rm d}V_a = \sum_i\int_M \lambda_i\,{\rm d}V_a = 0$ when $\bar Z=0$.
\end{remark}

Next theorem generalizes integral formula \eqref{eq61b-init}, using approach of foliated Randers spaces:
that is given a Riemannian space $(M, a)$ with a vector field~$\beta^\sharp$ of small norm,
we associate a Randers space $(M,\alpha+\beta)$.
Recall that $F=\alpha+\beta$ is Berwald if and only if $\bar\nabla\beta^\sharp=0$.
In this case, the Finsler metric and the source metric $\langle\cdot,\cdot\rangle$ have equal Riemann curvatures:
$R_y=\bar R_y$ for $y\in TM_0$, see Remark~\ref{C-3beta}.

\begin{thm}\label{Cor-k1}
Let a Riemannian manifold $(M, a)$ admits a non-trivial parallel vector field~$\beta^\sharp$
(say, $\|\beta^\sharp\|_\alpha<1$), which is nowhere orthogonal to a codimension-one foliation $\calf$. Then
\begin{eqnarray}\label{Eq-sigma_gen}
\nonumber
 &&\hskip-10mm \int_M c^{-2}\Big( \sigma_2(\bar A+c\,C^\sharp_\nu)
 +\big(\frac{c-2\,\hat c}{c\,\hat c}\,\<\bar A(\beta^{\sharp\top}),\beta^{\sharp}\>
 -\frac{\hat c-c}{c^2\hat c}\,\beta(\bar Z)\big)\,\sigma_1(\bar A+c\,C^\sharp_\nu) \\
\nonumber
 &&\hskip-10mm
 +\,\frac{\hat c-c}{c\,\hat c\,(1-c^2)}\,\<(\bar A+c\,C^\sharp_\nu)(\beta^{\sharp\top}), \beta^{\sharp\top}\>
  \<\bar A(\beta^{\sharp\top}), \beta^{\sharp\top}\>
 -\frac{c(c-2\,\hat c)^{\,2}(1-c^2)}{4\,c\,\hat c^{\,2}}\,\|\bar Z^{\bot\beta}\|^2_\alpha \\
\nonumber
 && \hskip-10mm
 +\frac{c-2\,\hat c}{c\,\hat c\,(1-c^{2})}\,\beta(\bar Z)\,\<(\bar A+c\,C^\sharp_\nu)(\beta^{\sharp\top}),\beta^{\sharp\top}\>
 -\frac{1-(c-2\,\hat c)^{\,2}}{4\,\hat c^{\,2}}\,\|\bar A(\beta^{\sharp\top})^{\bot\beta}\|_\alpha^2  \\
\nonumber
 &&\hskip-10mm
 -\,\frac{(c-2\,\hat c)(1-c^2+2\,c\,\hat c)}{2\,\hat c^{\,2}}\,\<\bar A(\beta^{\sharp\top})^{\bot\beta},\bar Z^{\bot\beta}\>
 -\frac{1+c^2-2\,c\,\hat c}{2\,\hat c}\,\<\bar A(\beta^{\sharp\top})^{\bot\beta}, C^\sharp_\nu(\beta^{\sharp\top})^{\bot\beta}\>\\
 &&\hskip-10mm -\,\frac{(c-2\,\hat c)(1+c^2)}{2\,\hat c}\,\<C^\sharp_\nu(\beta^{\sharp\top})^{\bot\beta},\bar Z^{\bot\beta}\>
 -\frac12\,\overline{\Ric}_{\,N}\Big)\,{\rm d}V_a = 0\,.
\end{eqnarray}
 Furthermore, if $\,\beta(N)=\const$, $N$ being a unit normal to $\calf$, then \eqref{Eq-sigma_gen} reads
\begin{eqnarray}\label{Eq-sigma_gen-const}
\nonumber
 &&\hskip-10mm \int_M \Big( c\tr(C^\sharp_\nu)\,\sigma_1(\bar A) -c\tr(\bar A C^\sharp_\nu)
 -\frac{1-(c-2\,\hat c)^{\,2}}{4\,\hat c^{\,2}}\,\|\bar A(\beta^{\sharp\top})\|_\alpha^2 \\
\nonumber
 &&\hskip-10mm
 -\,\frac{(c-2\,\hat c)(1-c^2+2\,c\,\hat c)}{2\,\hat c^{\,2}}\,\<\bar A(\beta^{\sharp\top}),\bar Z\>
 -\frac{1+c^2-2\,c\,\hat c}{2\,\hat c}\,\<\bar A(\beta^{\sharp\top}), C^\sharp_\nu(\beta^{\sharp\top})\>\\
 &&\hskip-10mm
 -\,\frac{c(c-2\,\hat c)^{\,2}(1-c^2)}{4\,c\,\hat c^{\,2}}\,\|\bar Z\,\|^2_\alpha
 -\frac{(c-2\,\hat c)(1+c^2)}{2\,\hat c}\,\<C^\sharp_\nu(\beta^{\sharp\top}),\bar Z\>
 \Big)\,{\rm d}V_a = 0\,.
\end{eqnarray}
\end{thm}

\proof
Note that $c<1$ when $\beta^{\sharp\top}\!\ne0$ on a Randers space $(M,\alpha+\beta)$.
For $\bar\nabla\beta^\sharp=0$ we get $(\bar\nabla_n\,\beta^{\sharp\top})^\top=-\beta(N)(\bar A(\beta^{\sharp\top})+c\bar Z)$.
 By~\eqref{E-dt-A2}$_1$ and \eqref{E-A-bar-A}
with $U=\frac{c(c-2\,\hat c)}{\hat c}\bar Z-\frac{\hat c-c}{\hat c}\,\bar A(\beta^{\sharp\top})$,
we have $A=A^g+C^\sharp_\nu$ and $cA^g=\bar A+A_1+A_2+A_3$, where
\begin{eqnarray*}
 A_1=U_1^\flat\otimes\beta^{\sharp\top},\quad
 A_2=U_2\otimes\beta^\top,\quad
 A_3= a_3 \,\beta^{\top}\!\otimes\beta^{\sharp\top}
\end{eqnarray*}
are rank 1 matrices,
 $a_3= \frac{\beta(U)}{c^2(1-c^2)}
 =\frac{c-2\,\hat c}{c\,\hat c(1-c^2)}\,\beta(\bar Z) -\frac{\hat c-c}{c^2\hat c(1-c^2)}\,\<\bar A(\beta^{\sharp\top}),\beta^{\sharp}\>$
and
\[
 U_1 = \frac1{2\,c\,\hat c}\,\big(\bar A(\beta^{\sharp\top}) +(c-2\,\hat c)\bar Z\big)^{\bot\beta},\quad
 U_2 = \frac{c-2\,\hat c}{2\,\hat c}\,\big(\bar A(\beta^{\sharp\top}) +c\,\bar Z\big)^{\bot\beta}.
\]
Thus, $c\,A=\bar A+c\,C^\sharp_\nu+A_1+A_2+A_3$. We have $\sigma_1(A_1)=\sigma_1(A_2)=\sigma_2(A_i)=0$~and
\begin{eqnarray*}
 \tr(A_1 A_2) \eq \<U_1,U_2\>\,\beta(\beta^{\sharp\top})
 =\frac{(c-2\,\hat c)(1-c^2)}{4\,c\,\hat c^{\,2}}\,\|\bar A(\beta^{\sharp\top})^{\bot\beta}\|^2_\alpha \\
 \minus \frac{(\hat c -c)(c-2\,\hat c)(1-c^2)}{2\,c\,\hat c^{\,2}}\,\<\bar A(\beta^{\sharp\top})^{\bot\beta},\bar Z\>
 +\frac{c(c-2\,\hat c)^2(1-c^2)}{4\,c\,\hat c^{\,2}}\,\|\bar Z^{\bot\beta}\|^2_\alpha,\\
 \tr(A_1(\bar A+c\,C^\sharp_\nu)) \eq \<U_1,\,(\bar A+c\,C^\sharp_\nu)(\beta^{\sharp\top})\>
  = \frac{1}{2\,c\,\hat c}\,\|\bar A(\beta^{\sharp\top})^{\bot\beta}\|^2_\alpha
  + \frac{1}{2\,\hat c}\,\<\bar A(\beta^{\sharp\top})^{\bot\beta},\,C^\sharp_\nu(\beta^{\sharp\top})\>\\
 \plus \frac{c-2\,\hat c}{2\,c\,\hat c}\,\<\bar A(\beta^{\sharp\top})^{\bot\beta},\,\bar Z\>
 +\frac{c-2\,\hat c}{2\,\hat c}\,\<\bar Z^{\bot\beta},\,C^\sharp_\nu(\beta^{\sharp\top})\> ,\\
 \tr(A_2(\bar A+c\,C^\sharp_\nu)) \eq
 \<U_2,\,(\bar A+c\,C^\sharp_\nu)(\beta^{\sharp\top})\>
 = \frac{c-2\,\hat c}{2\,\hat c}\,\|\bar A(\beta^{\sharp\top})^{\bot\beta}\|^2_\alpha
 +\frac{c(c-2\,\hat c)}{2\,\hat c}\,\<\bar A(\beta^{\sharp\top})^{\bot\beta},\,\bar Z\>\\
 \plus \frac{c(c-2\,\hat c)}{2\,\hat c}\,\<\bar A(\beta^{\sharp\top})^{\bot\beta},\,C^\sharp_\nu(\beta^{\sharp\top})\>
 +\frac{c^2(c-2\,\hat c)}{2\,\hat c}\,\<\bar Z^{\bot\beta},\,C^\sharp_\nu(\beta^{\sharp\top})\> ,\\
 \tr(A_3(\bar A+c\,C^\sharp_\nu)) \eq
 \frac{c-2\,\hat c}{c\,\hat c\,(1-c^{2})}\,\beta(\bar Z)\<(\bar A+c\,C^\sharp_\nu)(\beta^{\sharp\top}),\beta^{\sharp\top}\> \\
 \minus \frac{\hat c-c}{c\,\hat c\,(1-c^2)}\,\<\bar A(\beta^{\sharp\top}),\beta^{\sharp}\>
 \<(\bar A+c\,C^\sharp_\nu)(\beta^{\sharp\top}),\,\beta^{\sharp\top}\>,\\
 \tr(A_3 A_1) \eq a_3\<U_1,\,\beta^\sharp\>\,\beta(\beta^{\sharp\top}) = 0 ,\qquad
   \tr(A_3 A_2) = a_3\<U_2,\,\beta^\sharp\>\,\beta(\beta^{\sharp\top}) = 0 .
\end{eqnarray*}
Recall the following identity for square matrices:
\begin{eqnarray*}
 \sigma_2(\sum\nolimits_{\,i} A_i) \eq\sum\nolimits_{\,i}\sigma_2(A_i)+\sum\nolimits_{\,i<j} \big((\tr A_i)(\tr A_j) -\tr(A_i A_j)\big).
\end{eqnarray*}
By the above we obtain
\begin{eqnarray*}
  c^2\sigma_2(A) \eq \sigma_2(\bar A+c\,C^\sharp_\nu) +\sigma_1(A_3)\,\sigma_1(\bar A+c\,C^\sharp_\nu)\\
 \minus \tr(A_1 A_2+A_1(\bar A+c\,C^\sharp_\nu)+A_2(\bar A+c\,C^\sharp_\nu) +A_3(\bar A+c\,C^\sharp_\nu)),
\end{eqnarray*}
where $\sigma_1(A_3)=a_3(1-c^2)$ and
\begin{eqnarray*}
 &&\tr(A_1 A_2+A_1(\bar A+c\,C^\sharp_\nu)+A_2(\bar A+c\,C^\sharp_\nu) +A_3(\bar A+c\,C^\sharp_\nu) ) \\
 && = -\frac{\hat c-c}{c\,\hat c\,(1-c^2)}\,\<(\bar A+c\,C^\sharp_\nu)(\beta^{\sharp\top}), \beta^{\sharp\top}\>
  \<\bar A(\beta^{\sharp\top}), \beta^{\sharp\top}\> \\
 && -\,\frac{c-2\,\hat c}{c\,\hat c\,(1-c^{2})}\,\beta(\bar Z)\,\<(\bar A+c\,C^\sharp_\nu)(\beta^{\sharp\top}),\beta^{\sharp\top}\>
 +\frac{c(c-2\,\hat c)^2(1-c^2)}{4\,c\,\hat c^{\,2}}\,\|\bar Z^{\bot\beta}\|^2_\alpha\\
 && +\,\frac{1-(c-2\,\hat c)^{\,2}}{4\,\hat c^{\,2}}\,\|\bar A(\beta^{\sharp\top})^{\bot\beta}\|_\alpha^2
  +\frac{(c-2\,\hat c)(1-c^2+2\,c\,\hat c)}{2\,\hat c^{\,2}}
  \,\<\bar A(\beta^{\sharp\top})^{\bot\beta}, \bar Z\>\\
 && +\,\frac{1+c^2-2\,c\,\hat c}{2\,\hat c}\,\<\bar A(\beta^{\sharp\top})^{\bot\beta}, C^\sharp_\nu(\beta^{\sharp\top})\>
 +\frac{(c-2\,\hat c)(1+c^2)}{2\,\hat c}\,\<\bar Z^{\bot\beta}, C^\sharp_\nu(\beta^{\sharp\top})\> .
\end{eqnarray*}
The condition $\bar\nabla\beta^\sharp=0$ implies $\|\beta^\sharp\|_\alpha=\const$ and
$\bar R(X,Y)\beta^\sharp=0\ (X,Y\in TM)$. Using equality
\[
 \overline{\Ric}_{\,n}=
 \overline{\Ric}_{\,\hat c\,N-\beta^{\sharp}}=\hat c^{\,2}\,\overline{\Ric}_{\,n} +\overline{\Ric}_{\,\beta^\sharp}
 -2\,\hat c\sum\nolimits_{\,i}\bar R(N,b_i,\beta^\sharp,b_i),
\]
we obtain $\overline{\Ric}_{\,\nu}=(c\,\hat c)^{\,-2}\,\overline{\Ric}_{\,N} =c^{-2}\,\overline{\Ric}_{\,N}$.
 From \eqref{eq61}$_2$ for $F$, where the volume form is
${\rm d}V_F=(1-\|\beta^{\sharp}\|^{2}_\alpha)^{\frac{m+2}2}\,{\rm d}V_a$, see \eqref{E-F001vol},
we find \eqref{Eq-sigma_gen}. Since
$\lim\limits_{\,\beta\to0}\bar A(\beta^\sharp)^{\bot\beta}=0$,
\eqref{Eq-sigma_gen} reduces to \eqref{eq61b-init} when $\beta\to0$.

 If $\,\beta(N)=\const$ then $\beta(\bar Z)=0$ and $\<\bar A(\beta^{\sharp\top}),\beta^{\sharp}\>=0$:
\begin{eqnarray}\label{E-Abb-const}
\nonumber
 && 0=\<\bar\nabla_N\,\beta^\sharp,\,N\>
 =\<\bar\nabla_N\,(\beta^{\sharp\top}+\beta(N)N),\,N\>
 =-\<\beta^{\sharp\top},\,\bar\nabla_N N\> =-\<\beta^\sharp,\,\bar Z\>,\\
 && 0=\<\bar\nabla_{\beta^{\sharp\top}}\beta^\sharp,\,N\>
 =\<\bar\nabla_{\beta^{\sharp\top}}(\beta^{\sharp\top}+\beta(N)N),\,N\>
 =-\<\bar A(\beta^\sharp),\beta^\sharp\>.
\end{eqnarray}
Hence, and by Lemma~\ref{L-AplusB} for $\sigma_2(\bar A+c\,C^\sharp_\nu)$
and by \eqref{eq61b-init}, we reduce \eqref{Eq-sigma_gen} to \eqref{Eq-sigma_gen-const}.
\qed

\smallskip

Remark that a parallel vector field $\beta^\sharp$ forms a constant angle with (the leaves of) $\calf$
if and only if $\beta(N)=\const$ (e.g. $\beta(N)=0$) and $\|\beta^{\sharp\top}\|_\alpha=\const$.

\begin{corollary}\label{C2-gen}
Assume that a Riemannian manifold $(M, a)$ admits a non-trivial parallel vector field~$\beta^\sharp$,
which forms a constant angle with the leaves of a Riemannian $(\bar Z = 0)$ foliation $\calf$,
and $2\beta(N)+c\ne1$.
Then $\bar A(\beta^{\sharp\top})=0$ at any point of~$M$.
If, in addition, $\calf$ is totally umbilical $(\bar A = \bar H\cdot I_m)$ then $\calf$ is totally geodesic.
\end{corollary}

\proof
Let $\|\beta^\sharp\|_\alpha<1$.
By conditions and Corollary~\ref{C-sharp-n}(i), $\bar Z=0$ yields $C^\sharp_\nu=0$ on a Randers space $(M,\alpha+\beta)$.
Since $c$ and $\hat c$ are constant, \eqref{Eq-sigma_gen-const} reads
\begin{equation*}
 \frac{1-(2\,\hat c-c)^{\,2}}{4\,\hat c^{\,2}}\int_M \|\bar A(\beta^{\sharp\top})\|_\alpha^2\,{\rm d}V_a = 0.
\end{equation*}
By conditions, the factor $1-(2\,\hat c-c)^2$ is nonzero.
This yields $\bar A(\beta^{\sharp\top})=0$ on $M$.
If $\calf$ is a totally umbilical foliation then $0=\<\bar A(\beta^{\sharp\top}),\beta^{\sharp}\>=\bar H\|\beta^{\sharp\top}\|^2_\alpha$, hence $\bar H=0$.
\qed

\begin{corollary}\label{C2-k1}
Let a Riemannian manifold $(M, a)$
with a codimension-one foliation $\calf$ admits a nonzero parallel vector field~$\beta^\sharp\in\Gamma(T\calf)$
(say, $\|\beta^\sharp\|_\alpha<1$). Then
\begin{eqnarray}\label{Eq-sigma_gen-b}
\nonumber
 &&\hskip-10mm \int_M \Big( c\tr(C^\sharp_\nu)\,\sigma_1(\bar A) -c\tr(\bar A C^\sharp_\nu)
 -\frac{1-c^2}{4\,c^{2}}\,\|\bar A(\beta^{\sharp})\|_\alpha^2
 +\frac{1+c^2}{2\,c}\,\<\bar A(\beta^{\sharp}), \bar Z\> \\
 &&\hskip-10mm -\frac{1-c^2}{4}\,\|\bar Z\|^2_\alpha
 {-}\frac{1-c^2}{2\,c}\,\<\bar A(\beta^{\sharp}), C^\sharp_\nu(\beta^{\sharp})\>
 {+}\frac{1+c^2}{2}\,\<C^\sharp_\nu(\beta^{\sharp}),\,\bar Z\> \Big)\,{\rm d}V_a = 0.
\end{eqnarray}
If, in addition, $\bar Z=0$ then $\bar A(\beta^\sharp)=0$ at any point of $M$.
\end{corollary}

\proof By conditions, $c=\const<1$ on a Randers space $(M,\alpha+\beta)$.
Since $\bar\nabla\beta^\sharp=0$, we obtain $\beta(\bar Z)=0$ and $\<\bar A(\beta^\sharp),\beta^\sharp\>=0$:
\begin{eqnarray*}
 && 0=\<\bar\nabla_N\,\beta^\sharp,\,N\>=-\<\beta^\sharp,\,\bar\nabla_N N\> =-\<\beta^\sharp,\,\bar Z\>,\\
 && 0=\<\bar\nabla_{\beta^\sharp}\,\beta^\sharp,\,N\>=-\<\beta^\sharp,\,\bar\nabla_{\beta^\sharp}\,N\>
 =-\<\bar A(\beta^\sharp),\beta^\sharp\>.
\end{eqnarray*}
Comparing \eqref{Eq-sigma_gen} with \eqref{eq61b-init} for Riemannian metric $\langle\cdot,\cdot\rangle$
and $\beta(N)=0$ and applying Lemma~\ref{L-AplusB} to $\sigma_2(\bar A+c\,C^\sharp_\nu)$,
we reduce it to \eqref{Eq-sigma_gen-b}.
For $\bar Z=0$, by Corollary~\ref{C-sharp-n}(i) with $\beta(N)=0$, we have $C^\sharp_\nu=0$; hence, \eqref{Eq-sigma_gen-b} reads
$\frac{1-c^2}{4\,c^{2}}\int_M \|\bar A(\beta^{\sharp})\|_\alpha^2\,{\rm d}V_a = 0$.
Since $c<1$, this yields $\bar A(\beta^{\sharp})=0$ on $M$.
\qed

\begin{remark}\rm
Using formula in \cite[Lemma~4.2.2]{cs} for $\Ric-\overline{\Ric}$, see also Remark~\ref{R-22},
one may generalize \eqref{Eq-sigma_gen} (completing it with more terms) for
Randers spaces without additional condition $\bar\nabla\beta^\sharp=0$.
\end{remark}

\subsection{Around Brito-Langevin-Rosenberg formula}

Results of this section are valid for a codimension-one foliation and 1-form with singularities:
in Theorem~\ref{T-Randers-Kconst} and Corollary~\ref{C-Randers-A0} (according to \cite[Theorem~2 and Corollary~4]{rw2}
and Lemma~\ref{L-LuW-2}).

Recall that the \textit{Newton transformations} $T_r(A)\ (0\le r\le m)$ of an $m\times m$ matrix $A$ (see \cite{rw1}) are defined either inductively by
$T_0(A)=\,I_m$ and $T_r(A)=\sigma_r(A)\,I_m-A\,T_{r-1}(A)\ (r\ge1)$.
Note that $T_r(\lambda\,A)=\lambda^r\,T_r(A)$ for $\lambda\ne0$
and $\tr(T_r(A))=(m-r)\,\sigma_r(A)$.
Observe that if~a rank-one matrix $A:=U\otimes\beta$ (and similarly for $A:=\omega\otimes\beta^\sharp$)
has  trace zero, i.e., $\beta(U)=0$, then we~have
\[
 A^2=U(\beta^\sharp)^t \cdot U(\beta^\sharp)^t = U\beta(U)\,(\beta^\sharp)^t=\beta(U)\,A=0.
\]
Define the quantity
\[
 \delta:=-\frac12\,c^{-1}\hat c^{\,-2}(\hat c\,N-\beta^{\sharp})(c\,\hat c).
\]

In this section we assume that our Randers space is Berwald, and $\beta^\sharp$ is nowhere orthogonal to $\calf$
(for $\beta$ orthogonal to $\calf$ see Remark~\ref{sec:beta0}):
\begin{equation}\label{E-nabla-beta}
 \bar\nabla\beta^\sharp=0,\qquad  \beta^{\sharp\top}\ne0.
\end{equation}
If, in addition, $\langle\cdot,\cdot\rangle$ has constant curvature $\bar K$ then
$\bar K=0$
(because only flat space forms admit parallel vector fields).
Indeed, since $\bar R(x,y)z=\bar K(\,\<y,z\>\,x-\<x,z\>\,y\,)$, on $T\calf$ we have
\[
 \bar R_N=\bar K \,I_m,\quad
 \bar R_{\beta^\sharp}=(1-c^2)\bar K \,I_m,\quad
 \bar R(\cdot,N)\beta^\sharp=0.
\]
If $\,\bar\nabla\beta=0$ then $\bar R(U,\beta^\sharp,\beta^\sharp,U) = 0$ and $\bar K(U\wedge\beta^\sharp) = 0$
for all nonzero $U\perp\beta^\sharp$; hence, $\bar K=0$.

\begin{thm}\label{T-Randers-Kconst}
Let $(M^{m+1},\alpha+\beta)$ be a codimension-one foliated closed Ran\-ders-Berwald space
with conditions \eqref{E-nabla-beta} and constant sectional curvature $\bar K=0$ of $\langle\cdot,\cdot\rangle$.
Then
\begin{eqnarray}\label{E-IF-Randers-k}
\nonumber
 &&\hskip-8mm \int_M \Big( \delta\,(m-k+1)\,\sigma_{k-1}(\bar A)
 +\sum\nolimits_{j>0}\sigma_{k-j,j}(\bar A+\delta\,I_m,\ c\,C^\sharp_\nu)\\
\nonumber
 &&\hskip-6mm +\<T_{k-1}(\bar A+\delta\,I_m{+}c\,C^\sharp_\nu)(\beta^{\sharp\top}),\ U_1\>
 +\beta\big(T_{k-1}(\bar A+\delta\,I_m{+}c\,C^\sharp_\nu+U_1^\flat\otimes\beta^{\sharp\top})(U_2)\big)\\
 &&\hskip-6mm +\,a_3\,
 \beta\big(T_{k-1}(\bar A+\delta\,I_m +c\,C^\sharp_\nu+U_1^\flat\otimes\beta^{\sharp\top}
 \!+U_2\otimes\beta^\top)(\beta^{\sharp\top})\big)\Big){\rm d}V_a = 0,
\end{eqnarray}
where $1\le k\le m$,
$a_3=\frac{c-2\,\hat c}{c\,\hat c(1-c^2)}\,\beta(\bar Z) -\frac{\hat c-c}{c^2\hat c(1-c^2)}\,\<\bar A(\beta^{\sharp\top}),\beta^{\sharp}\>$
and
\[
 U_1 = \frac1{2\,c\,\hat c}\,\big(\bar A(\beta^{\sharp\top}) +(c-2\,\hat c)\bar Z\big)^{\bot\beta},\quad
 U_2 = \frac{c-2\,\hat c}{2\,\hat c}\,\big(\bar A(\beta^{\sharp\top}) +c\,\bar Z\big)^{\bot\beta}.
\]
Furthermore, if $\,\beta(N)=\const$ and $\bar Z=0$ then
\begin{equation*}
 \int_M \frac{1+c^2-2\,c\,\hat c}{2\,c\,\hat c}\,
 \big\<T_{k-1}(\bar A)(\beta^{\sharp\top}),\ \bar A(\beta^{\sharp\top})^{\bot\beta}\big\>\,{\rm d}V_a\!=0.
\end{equation*}
\end{thm}

\proof As was shown, $\bar K=0$, and $R_y = \bar R_y=0$ for $y\in TM_0$.
By assumptions, $c<1$ and $\|\beta\|_\alpha=\const$.
By \eqref{E-dt-A2} and \eqref{E-A-bar-A},
\[
 c\,A = c\,A^g+c\,C^\sharp_\nu = \bar A+\delta\,I_m +c\,C^\sharp_\nu+A_1+A_2 +A_3,
\]
where $A_i$ are three rank $\le 1$ matrices,
\[
 A_1=U_1^\flat\otimes\beta^{\sharp\top},\ \
 A_2=U_2\otimes\beta^\top,\ \
 A_3=a_3\,\beta^\top\!\otimes\beta^{\sharp\top}.
\]
By Lemma~\ref{L-AplusB} (with $C=\bar A+\delta\,I_m$ and $D=c\,C^\sharp_\nu$) of Appendix, we have
\begin{eqnarray}\label{E-2sigma-k}
\nonumber
 && c^k\sigma_k(A) = \sigma_k(\bar A+\delta\,I_m) +\sum\nolimits_{j>0}\sigma_{k-j,j}(\bar A+\delta\,I_m, c\,C^\sharp_\nu)  \\
\nonumber
 && +\,U_1^\flat(T_{k-1}(\bar A+\delta\,I_m+c\,C^\sharp_\nu)(\beta^{\sharp\top}))
 +\beta(T_{k-1}(\bar A+\delta\,I_m+c\,C^\sharp_\nu+A_1)(U_2))\\
 && +\,a_3\,\beta(T_{k-1}(\bar A+\delta\,I_m+c\,C^\sharp_\nu{+}A_1{+}A_2)(\beta^{\sharp\top})).
\end{eqnarray}
Recall that  ${\rm d}V_F=(1-\|\beta^{\sharp}\|^{2}_\alpha)^{\frac{m+2}2}{\rm d}V_a$, see \eqref{E-F001vol}$_1$.
Comparing \eqref{eq5f-b} when $K=0$ with $\int_M \sigma_k(\bar A_p)\,{\rm d}V_a = 0$,
and using $\sigma_k(\bar A+\delta\,I_m)=\sigma_k(\bar A) +\delta(m-k+1)\,\sigma_{k-1}(\bar A)$,
see Lemma~\ref{L-AplusB}, and \eqref{E-2sigma-k}, we find \eqref{E-IF-Randers-k}.

If $\,\beta(N)=\const$ and $\bar Z=0$ then $c$ and $\hat c$ are constant; hence, $\delta=0$
and $\<\bar A(\beta^{\sharp\top}),\,\beta^{\sharp\top}\>=0$ and $a_3=0$.
Thus, second claim follows from Corollary~\ref{C-sharp-n}(i) and \eqref{E-IF-Randers-k}.
\qed

\begin{example}\label{Ex-k1}\rm
Let conditions of Theorem~\ref{T-Randers-Kconst} hold. For $k=1$, \eqref{E-IF-Randers-k} yields the Reeb type formula
\begin{equation*}
 \int_M\big(\,c\tr(C^\sharp_\nu)+m\,\delta +\frac{c-2\,\hat c}{c\,\hat c}\,\beta(\bar Z)
 -\frac{\hat c -c}{c^2\hat c}\,\<\bar A(\beta^{\sharp\top}),\beta^{\sharp}\>\big)\,{\rm d}V_a = 0
\end{equation*}
(since $\beta(U_1)=\beta(U_2)=0$),
which for $\bar Z=0$ reads $\int_M\frac{\hat c -c}{c^2\hat c}\,\<\bar A(\beta^{\sharp\top}),\beta^{\sharp}\>\,{\rm d}V_a=0$,
see also \eqref{E-IF1-Randers-gen}.
\end{example}

\begin{corollary}\label{C-Randers-A0}
Assume that  $(M^{m+1},\alpha+\beta)$ is a closed Randers-Berwald space of constant sectional curvature~$\bar K=0$ of $\langle\cdot,\cdot\rangle$ endowed with a codimension-one totally geodesic $($for our Riemannian metric $a)$ foliation and conditions \eqref{E-nabla-beta} hold. Then for $1\le k\le m$ $($for $k=1$, see also Example~\ref{Ex-k1}$)$
\begin{eqnarray}\label{E-IF-Randers-k-tg}
\nonumber
 &&\hskip-8mm\int_M\!\Big( c^k\sigma_{k}(C^\sharp_\nu)
 +\frac{c-2\,\hat c}{2\,c\,\hat c}\,\<T_{k-1}(C^\sharp_\nu+\delta\,I_m)(\beta^{\sharp\top}),\bar Z^{\bot\beta}\> \\
\nonumber
 &&\hskip-8mm +\,\frac{c(c-2\,\hat c)}{2\,\hat c}\,\beta\big(T_{k-1}(c\,C^\sharp_\nu+\delta\,I_m
 +\frac{c-2\,\hat c}{2\,c\,\hat c}\,(\bar Z^{\bot\beta})^\flat\otimes\beta^{\sharp\top})(\bar Z^{\bot\beta})\big)\\
\nonumber
 &&\hskip-8mm +\,\frac{c-2\,\hat c}{c\,\hat c\,(1-c^2)}\,\beta(\bar Z)\,\beta\big(T_{k-1}(c\,C^\sharp_\nu+\delta\,I_m
 +\frac{c-2\,\hat c}{2\,c\,\hat c}\,(\bar Z^{\bot\beta})^\flat\otimes\beta^{\sharp\top}\\
 &&\hskip-8mm +\frac{c(c-2\,\hat c)}{2\,\hat c}\,\bar Z^{\bot\beta}\otimes\beta^\top)(\beta^{\sharp\top})\big)\Big){\rm d}V_a=0.
\end{eqnarray}
\end{corollary}

\proof This follows from \eqref{E-IF-Randers-k} with $\bar A=0$ and
\[
 U_1=\frac{c-2\,\hat c}{2\,c\,\hat c}\,\bar Z^{\bot\beta},\quad
 U_2=\frac{c(c-2\,\hat c)}{2\,\hat c}\,\bar Z^{\bot\beta},\quad
 a_3=\frac{c-2\,\hat c}{c\,\hat c\,(1-c^2)}\,\beta(\bar Z).
\]
For~$\beta(N)=0$, \eqref{E-IF-Randers-k-tg} reduces to formula (4.23) in \cite{rw3}.
\qed

\smallskip

For~$\beta(N)=0$, \eqref{E-IF-Randers-k-tg} reduces to formula (4.23) in \cite{rw3}.
Similar integral formulae exist for totally umbilical foliations (for $\beta(N)\equiv0$ see \cite[Corollary~4.7]{rw3}),
i.e. $\bar A=\bar H I_m$, where $\bar H=\frac1m\,\sigma_1(\bar A)$,~and
\[
 c\,A = c\,A^g+c\,C^\sharp_\nu = (\bar H+\delta)\,I_m +c\,C^\sharp_\nu+A_1+A_2 +A_3.
\]
Note that non-flat {closed} Riemannian manifolds of constant curvature do not admit such foliations.

\begin{remark}\label{sec:beta0}\rm
Let $\beta^\sharp=fN$ for a smooth function $f: M\to(-1,1)$. Then $c=1$ and $\beta(N)=f$.

1) By Lemma~\ref{L-c-value}, $n=N$, $\nu = \hat c^{\,-1}N$, and
\begin{eqnarray*}
 g(n,n) = \hat c^{\,2},\quad g(u,v) = \hat c\,\<u,v\>\quad(u,v\in T\calf),
\end{eqnarray*}
where $\hat c=1 + f$.
By \eqref{E-c-value0}, for arbitrary $u,v\in TM$ we have
\begin{equation*}
 g(u,v) = (1+f)\,(\,\<u,v\> +f\<N,u\>\,\<N,v\>\,).
\end{equation*}
Note that $\overline{\Div}\,^\top\bar Z=\overline{\Div}\,\bar Z+\<\bar Z,\bar Z\>$ and, see \eqref{E-tr-Def}:
 $\tr(\overline{\rm Def}_{fN})_{\,|\,T\calf}^\top =\<\bar\nabla_i(fN),b_i\> = - f\sigma_1(\bar A)$.
From \eqref{E-A-bar-A-initial} and Propositions~\ref{L-Dx} and \ref{P-ZZ} we obtain the following:
\begin{eqnarray}
\nonumber
 && A^g = \bar A -\frac12\,\hat c^{\,-2}N(f)\,I_m
 +\hat c^{\,-1}\,(\overline{\rm Def}_{fN})_{\,|\,T\calf}^\top,\\
\label{E-A-bar-A0-2}
 &&\sigma_1(A^g) = \hat c^{\,-1}\sigma_1(\bar A) -\frac m2\,\hat c^{\,-2} N(f),\\
\nonumber
 && Z = \hat c^{\,-1}\bar Z -\hat c^{\,-2}\,\bar\nabla^\top f,\qquad
 2\,C^{\,\sharp}_n = \hat c^{\,-3}f N(f)\,I_m.
\end{eqnarray}
Moreover, if $f=\const$ $($hence, $\hat c =\const)$ then
\begin{eqnarray*}
 A^g \eq \bar A +f\,\hat c^{\,-1}\,(\overline{\rm Def}_{N})_{\,|\,T\calf}^\top,\quad
 \sigma_1(A^g) = \hat c^{\,-1}\sigma_1(\bar A),\quad
 Z = \hat c^{\,-1}\bar Z,\quad
 C^{\,\sharp}_n = 0.
\end{eqnarray*}
The canonical volume forms of Riemannian metrics $g$ and $\langle\cdot,\cdot\rangle$ satisfy, see \eqref{E-F001vol}$_2$,
\begin{equation*}
 {\rm d} V_g = (1+f)^{\frac{m+2}2}{\rm d} V_a.
\end{equation*}

2) Let $F = \alpha + \beta$, where $\beta = f\cdot N^\flat$, be the corresponding Randers structure on $M$ and $g$ be the Riemannian metric on $M$ given by \eqref{E-c-value0}.
Comparing Reeb integral formula for metrics $\langle\cdot,\cdot\rangle$ and $g$ (or, just from \eqref{E-IF1-Randers-gen0})
and applying \eqref{E-A-bar-A0-2}, we~get:
\begin{eqnarray*}
 \int_M \sigma_{1}(A^g)\,{\rm d}V_g
 \eq \int_M(1+f)^{\frac{m+2}2}\big(\hat c^{\,-1}\sigma_1(\bar A) -\frac m2\,\hat c^{\,-2}N(f)\,\big)\,{\rm d}V_a\\
 \eq \int_M\Big( (1+f)^{\frac{m}2}\sigma_1(\bar A) -\frac m2\,(1+f)^{\frac{m-2}2}N(f)\,\Big)\,{\rm d}V_a,
\end{eqnarray*}
which vanishes for any $f$, see \eqref{E-reeb-f}. Hence, Theorem~\ref{T-1-1} yields a trivial identity in this case.

3) Let $(M^{m+1},\alpha+\beta)$ be a codimension-one foliated closed Randers space with constant
sectional curvature $\bar K$ of $\langle\cdot,\cdot\rangle$, and conditions $\beta^{\sharp\top}=0$ and $\bar\nabla\beta^\sharp=0$.
Then $\bar K=0$, see Theorem~\ref{T-Randers-Kconst}.
 Note that $\bar\nabla\beta^\sharp=0$ means $f=\const$ and $\bar\nabla N=0$. Hence, $C^\sharp_\nu=0$.
By \eqref{E-IF-Randers-k} with
$U_1 = \frac{c-2\,\hat c}{2\,c\,\hat c}$ and $U_2 = \frac{c(c-2\,\hat c)}{2\,\hat c}$,
\begin{eqnarray*}
 \int_M \big(\sigma_{k}(\bar A) +(-2)^{1-k}f\,(1+f)^{1-k}\,\<\,T_{k-1}(\bar A)(\bar Z),\ N\>\big){\rm d}V_a = 0.
\end{eqnarray*}
Since \eqref{eq5f-b} (for $\langle\cdot,\cdot\rangle$) and $\<\,T_{k-1}(\bar A)(\bar Z),\ N\>=0$, the above is satisfied for any $f$ and any $k>0$.
Again, Theorem~\ref{T-Randers-Kconst} yields a trivial identity in this case.
\end{remark}

\section{Appendix: Invariants of a set of quadratic matrices}
\label{sec:inv}

Here, we collect the properties of the invariants
$\sigma_{\lambda} (A_1, \ldots , A_k)$ of real matrices $A_i$ that generalize
the elementary symmetric functions of {a single matrix}~$A$.
Let $S_k$ be the group of all permutations of $k$ elements.
 Given arbitrary quadratic $m\times m$ real matrices $A_1, \ldots  A_k$
and the unit matrix $I_{m}$, one can consider the determinant $\det(I_{m}+t_1A_1+\ldots+t_kA_k)$
and express it as a polynomial of real variables ${\bf t}=(t_1, \dots  t_k)$. Given
 $\lambda =(\lambda_1, \ldots  \lambda_k)$, a sequence of nonnegative integers with
$|\lambda| := \lambda_1 + \ldots + \lambda_k\le m$, we shall
denote by $\sigma_{\lambda} (A_1, \ldots , A_k)$ its coefficient at
 ${\bf t}^{\lambda}=t_1^{\lambda_1}\cdot\ldots t_k^{\lambda_k}$:
\begin{equation*}
 \det(I_{m}+t_1A_1+\ldots+t_kA_k)=\sum\nolimits_{\,|\lambda|\,\le m}
 \sigma_{\lambda}(A_1, \ldots  A_k)\,{\bf t}^{\lambda}.
\end{equation*}
Evidently, the quantities $\sigma_{\lambda}$ are invariants of conjugation by $GL(m)$-matrices:
 $\sigma_{\lambda}(A_1,\ldots A_k)=\sigma_{\lambda}(QA_1Q^{-1},\ldots QA_kQ^{-1})$
for all $A_i$'s, ${\lambda}$'s and nonsingular $m\times m$ matrices $Q$.
Certainly, $\sigma_i(A)$ (for a single matrix $A$)
coincides with the $i$-th elementary symmetric polynomial of the
eigenvalues $\{k_j\}$ of~$A$. All the invariants $\sigma_{\lambda}$
can be expressed in terms of the traces of the matrices involved and their products.
 In the next lemmas, we collect properties of these invariants.

\begin{lem}[see \cite{rw2,rw3}]\label{lem1}
 For any $\lambda=(\lambda_1,\ldots \lambda_k)$ and any $m\times m$ matrices $A_i, A$  and $B$ one has

 (I) $\sigma_{\lambda}(0, A_2, \dots  A_{k})=0$ if $\lambda_1>0$
 and
 $\sigma_{0,\hat{\lambda}}(A_1, \dots  A_{k}) =\sigma_{\hat{\lambda}}(A_2,\ldots A_{k})$,
 where $\hat{\lambda}=(\lambda_2,\ldots \lambda_k)$,

 (II) $\sigma_{\lambda} (A_{s(1)},\ldots  A_{s(k)})=\sigma_{\lambda \circ s}
(A_1,\ldots A_k)$,
 where $s\in S_k$ and $\lambda\circ s=(\lambda_{s(1)},\ldots \lambda_{s(k)})$,

 (III) $\sigma_{\lambda}(I_{m},A_2,\ldots  A_k)
 = \genfrac{(}{)}{0pt}{}{m-|\hat{\lambda}|}{\lambda_1}\,
 \sigma_{\hat{\lambda}} (A_2,\ldots  A_k)$,

 (IV) $\sigma_{\lambda_1,\lambda_2,\,\hat{\lambda}}(A,A,A_3,\ldots  A_k)
 =\genfrac{(}{)}{0pt}{}{\lambda_1+\lambda_2}{\lambda_1}\,
 \,{\sigma_{\lambda_1+\lambda_2,\hat{\lambda}}}(A,A_3,\ldots  A_k)$,

 (V) $\sigma_{1,\hat{\lambda}}(A+B,A_2,\dots  A_{k}) =\sigma_{1,\hat{\lambda}}(A,A_2,\ldots
 A_{k})+\sigma_{1,\hat{\lambda}}(B,A_2,\ldots A_{k})$ and

 \qquad
 $\sigma_{\lambda}(a\/A_1,A_2,\dots  A_{k})
 =a^{\lambda_1}\sigma_{\lambda}(A_1,A_2,\ldots A_{k})$ if $a\in\RR\setminus\{0\}$.
\end{lem}

\begin{lem}[\cite{rw2,rw3}]
 For arbitrary matrices $B$, $C$ and $k,l>0$ we have
\begin{equation*}
 \sigma_{k,l}(B,C)=\sigma_k(B)\,\sigma_l(C)
 -\sum\nolimits_{\,i=1}^{\,\min(k,l)}\sigma_{k-i,l-i,i}(B,C,BC).
\end{equation*}
\end{lem}

\begin{lem}[\cite{rw3}]\label{L-AplusB}
Let $C,D,A_i\ (i\le s)$ be $m\times m$ matrices, ${\rm rank}\,A_i=1$. Then
\begin{eqnarray*}
 && \sigma_k(C+D+A_1+\ldots+A_s) = \sigma_k(C) +\sum\nolimits_{j>0}\sigma_{k-j,j}(C, D) \\
 \plus\tr(T_{k-1}(C+D)A_1) +\ldots+\tr(T_{k-1}(C+D+A_1+\ldots+A_{s-1})A_s).
\end{eqnarray*}
In particular (when $D=0$ and $s=1$),
 $\sigma_k(C+A) = \sigma_k(C) +\tr(T_{k-1}(C)A)$.
\end{lem}


\begin{thebibliography}{999.}%

\bibitem{arw2014}
 Andrzejewski K., Rovenski V. and Walczak P. {Integral formulas in foliation theory}, 73--82,
in ``\textit{Geometry and its Applications}", Springer Proc. in Math. and Statistics, 72, Springer,~2014.

\bibitem{brs}
D. Bao, C. Robles and Z. Shen, {Zermelo navigation on Riemannian manifolds},
J. Differential Geometry, 66 (2004), 377--435.

\bibitem{blr}
 Brito F., Langevin R. and  Rosenberg H. {Int\'{e}grales de courbure sur des vari\'{e}t\'{e}s feuillet\'{e}es}.
J.~Diff. Geom. {16} (1981), 19--50.

\bibitem{bw}
Brito F. and Walczak P. {On the energy of unit vector fields with isolated singularities}, Ann. Polon. Math.  {LXXIII.3} (2000),~269--274.


\bibitem{cs}
 Cheng X. and  Shen Z. \textit{Finsler geometry. An approach via Randers spaces},
Springer, 2012.

\bibitem{cs2}
 Chern S.S. and Shen Z. \textit{Riemann-Finsler geometry}, World Scientific, 2005.

\bibitem{lw}
Langevin R., Walczak P. {Conformal geometry of foliations}, Geom. Dedic.  {132} (2008),~135--178.

\bibitem{lw2}
Lu\.{z}y\'{n}czyk M., Walczak P.
{New integral formulae for two complementary orthogonal distributions on Riemannian manifolds},
Ann. Glob. Anal. Geom. 48 (2015), 195--209.

\bibitem{no}
 Nora T. {Seconde forme fondamentale d'une application et d'un feuilletage}, Th{\'e}se, l'Univ. de Limoges, 1983.

\bibitem{ra}
 Randers G. {On an asymmetrical metric in the four-space of general relativity}, Phys. Rev. 59 (1941), 195--199.

\bibitem{re}
 Reeb G. {Sur la courbure moyenne des vari\'{e}t\'{e}s int\'{e}grales d'une \'{e}quation de Pfaff} $\omega=0$.
 C. R. Acad. Sci. Paris {\bf 231}, 101--102 (1950)


\bibitem{rw1}
 Rovenski V. and Walczak P. \textit{Topics in extrinsic geometry of codimension-one foliations}, Springer 2011.

\bibitem{rw2}
 Rovenski V. and Walczak P. {Integral formulae on foliated symmetric spaces}, Math. Ann. {352} (2012), 223--237.

\bibitem{rw3}
Rovenski V. and Walczak P. {Integral formulae for codimension-one foliated Finsler spaces},
Balkan J. of Geometry and Its Appl. 21, No. 1 (2016), 76--102 (see ArXiv:1602.00610).


\bibitem{sh1}
 Shen Z. \textit{On Finsler geometry of submanifolds}, Math. Ann. {311} (1998), 549--576.

\bibitem{sh2}
 Shen Z. \textit{Lectures on Finsler geometry}, World Scientific Publishers, 2001.

\end{thebibliography}
\end{document}